\documentclass[final,3p,times]{elsarticle} 

\usepackage{amssymb}
\usepackage{amsmath}
\usepackage{algorithmic}
\usepackage{algorithm}
\usepackage{txfonts}
\usepackage{afterpage}
\usepackage{tikz}
\usepackage{xcolor}  
\usetikzlibrary{patterns} 
\usepackage{tikz-3dplot}
\usepackage{tabularx}
\usepackage{float}
\usepackage[pdftex]{pict2e}
\usepackage{lipsum}
\usepackage{amsfonts}
\usepackage{graphicx}
\usepackage{epstopdf}
\usepackage{subfig}
\usepackage{hyperref}  
\usepackage{booktabs}
\usepackage{multirow}
\usepackage{makecell}
\usepackage[pagewise]{lineno}

\usepackage{tabularx}
\usepackage{tabularray}
\usepackage{comment}
\usepackage{cases}

\newcommand{\bx}{\boldsymbol{x}}
\newcommand{\bxi}{\boldsymbol{\xi}}

\newcommand{\pF}{\mathcal{F}}

\newcommand{\dsR}{\mathbb{R} }

\usepackage{amsopn}

\ifpdf
  \DeclareGraphicsExtensions{.eps,.pdf,.png,.jpg}
\else
  \DeclareGraphicsExtensions{.eps}
\fi

\makeatletter
\newenvironment{breakablealgorithm}
  {
   \begin{center}
     \refstepcounter{algorithm}
     \hrule height.8pt depth0pt \kern2pt
     \renewcommand{\caption}[2][\relax]{
       {\raggedright\textbf{\ALG@name~\thealgorithm} ##2\par}%
       \ifx\relax##1\relax 
         \addcontentsline{loa}{algorithm}{\protect\numberline{\thealgorithm}##2}%
       \else 
         \addcontentsline{loa}{algorithm}{\protect\numberline{\thealgorithm}##1}%
       \fi
       \kern2pt\hrule\kern2pt
     }
  }{
     \kern2pt\hrule\relax
   \end{center}
  }
\makeatother

\newcommand{\warn}[1]{{\color{red}#1}}

\journal{}

\begin{document}

\begin{frontmatter}

\title{An adaptive adjoint-oriented neural network for\\ solving parametric optimal control problems with singularities}

\author[1]{Zikang Yuan}

\author[2]{Guanjie Wang}

\author[1]{Qifeng Liao}


\address[1]{School of Information Science and Technology, ShanghaiTech University, Shanghai, 201210, China}
\address[2]{School of Statistics and Mathematics, Shanghai Lixin University of Accounting and Finance, Shanghai, 201209, China}

\begin{abstract}
In this work, we present an adaptive adjoint-oriented neural network (adaptive AONN) for solving parametric optimal control problems governed by partial differential equations. The proposed method integrates deep adaptive sampling techniques with the adjoint-oriented neural network (AONN) framework. It alleviates the limitations of AONN in handling low-regularity solutions and enhances the generalizability of deep adaptive sampling for surrogate modeling without labeled data ($\text{DAS}^2$). The effectiveness of the adaptive AONN is demonstrated through numerical examples involving singularities.
\end{abstract}


\begin{keyword}
partial differential equations \sep parametric optimal control \sep adaptive sampling \sep surrogate modeling \sep singularities
\MSC[2020] 49M41 \sep 68T07 \sep 62D05


\end{keyword}

\end{frontmatter}



\section{Introduction} \label{section_intro}
Optimal control problems (OCPs) constrained by partial differential equations (PDEs) have found widespread applications across scientific and engineering disciplines, ranging from fluid dynamics~\cite{mathew2007optimal} and aerodynamic design~\cite{jameson1988aerodynamic,sun2020surrogate} to uncertainty quantification~\cite{xiu2010numerical,xiu2016stochastic} and control of  energy systems~\cite{lund2009energy}. These problems are often characterized by the need of optimizing certain objective functionals subject to systems governed by partial differential equations (PDEs). In practice, the objective functionals and/or the governing PDEs inherently involve various physical or geometric parameter configurations, leading to parametric OCPs constrained by PDEs. Such parameters typically arise from material properties, boundary conditions, applied control constraints, and the computational domain~\cite{karcher2018certified,negri2015reduced,stadler2009elliptic,milani2008reduced}.

Classical methods for PDE-constrained OCPs are largely based on local descent methods, such as Newton iterations~\cite{kelley1987quasi,sternberg2010memory} and the direct-adjoint looping (DAL) method~\cite{lions1971optimal,cacace2024reliable}, typically require the numerical solution of the governing PDEs at each iteration step. This process is often computationally intensive, and the presence of parameters in parametric OCPs further increases the computational burden. Substantial computational savings can be achieved if surrogate models are constructed to efficiently predict solutions with errors below a desired tolerance. A widely used method for this purpose is the reduced basis method (RBM)~\cite{quarteroni2015reduced, negri2013reduced, boyaval2010reduced, elman2013reduced,wang2024reduced}, which can provide accurate and stable approximations when the parametric solutions lie in a low-dimensional manifold~\cite{cohen2023nonlinear}. However, constructing the RBM approximation space typically requires a substantial amount of high-fidelity data, obtained through expensive numerical simulations of the governing PDEs. Another challenge with RBM is that, for high-dimensional problems, the dimension of the RBM approximation space is usually large, which may significantly reduce the efficiency of the method. 

In recent years, deep learning methods for solving PDEs have attracted increasing attention, with Physics-Informed Neural Networks (PINNs)~\cite{raissi2018hidden, raissi2019physics} and Deep Operator Network (DeepONet)~\cite{lu2019deeponet} emerging as two notable classes of neural network models. Both PINNs and DeepONet have been explored in the context of PDE-constrained OCPs. For example, the potential of PINNs to solve PDE-constrained OCPs is investigated in~\cite{mowlavi2023optimal}, by integrating the cost functional into the loss function of standard PINNs. In addition, operator learning methods have been proposed to construct surrogate models for the PDE solution operator, which can be integrated into solving PDE-constrained OCPs~\cite{lye2021iterative,hwang2022solving}.

Various learning-based methods have been extended to parametric OCPs constrained by PDEs~\cite{onken2022neural,demo2023extended,verma2024neural,zhang2025enhanced}. In particular, the adjoint-oriented neural network (AONN)~\cite{yin2024aonn} employs three distinct neural networks as parametric surrogate models for the state, control and adjoint functions, enabling the computation of optimal solutions in an all-at-once manner. To reduce training difficulty and handle complex constraints, the classical DAL method~\cite{lions1971optimal} is introduced into the iterative training framework, thereby avoiding penalty terms arising from the Karush-Kuhn-Tucker (KKT) system.

Although the AONN method leverages the strengths of both the DAL approach and deep learning techniques, designing effective random sampling strategies is crucial for training the neural networks within the AONN framework, as the statistical errors introduced by these random samples can dominate the overall approximation error, particularly for problems with low regularity or high dimensionality. 

In this work, we propose an adaptive AONN method, which integrates deep adaptive sampling techniques with the AONN framework, for solving parametric OCPs constrained by PDEs. More precisely, we construct surrogate models for the state, control, and adjoint functions using AONN, while sampling the geometric space and parameter space together with deep adaptive sampling for surrogate modeling without labeled data ($\text{DAS}^2$)~\cite{wang2024deep}. {In adaptive AONN,} we adaptively update the surrogate models and their training set alternately, with the training set update achieved by approximating the probability distribution induced by the residual profile. At each stage, the solution obtained through AONN enhances the ability of $\text{DAS}^2$ to approximate the target distribution, while the samples generated by $\text{DAS}^2$ aid in training the surrogate models within AONN, creating a mutually reinforcing computational synergy. Once the surrogate models are constructed and the sampling network is trained, the corresponding optimal control function and state function can be efficiently obtained for any set of parameters without the need of repeatedly solving the optimal control system.

The remainder of this paper is organized as follows. In Section~\ref{section_setup}, we introduce the problem settings and the AONN method for solving PDE-constrained OCPs via surrogate modeling. The details of the adaptive sampling strategy are presented in Section~\ref{section_sampling}. Subsequently, we propose the adaptive AONN for solving PDE-constrained OCPs in Section~\ref{section_method}. Numerical experiments are provided in Section~\ref{section_experiments} to demonstrate the effectiveness of the proposed method. Finally, Section~\ref{section_conclude} concludes the paper.

\section{Parametric optimal control problems and AONN}\label{section_setup}
This section describes the parametric optimal control problems considered in this study and introduces the adjoint-oriented neural network for solving them.

\subsection{Problem setup}
Let $\bxi \in \Gamma \subset \mathbb{R}^d$ be a vector collecting $d$ parameters, and  $\Omega(\bxi) \subset \dsR^n$ be a bounded and connected spatial domain with the boundary~$\partial\Omega(\bxi)$ that depends on $\bxi$. Let $\bx \in \Omega(\bxi)$ denote a spatial variable. In this work, we consider the following parametric optimal control system:
\begin{equation}\label{OCP}
    \text{OCP}(\bxi): \quad \left\{\begin{aligned}
        &\qquad \ \min \limits_{(y(\bx, \bxi), u(\bx, \bxi)) \in Y \times U} J\big(y(\bx, \bxi), u(\bx, \bxi); \bxi\big), \\ 
        & \text{subject to} \quad \pF\big(y(\bx, \bxi), u(\bx, \bxi), \bxi \big) = 0 \ \text{in} \  \Omega(\bxi), \ \text{and}\ u(\bx, \bxi) \in U_{ad}(\bxi), \\
        \end{aligned}\right. \\
\end{equation}
where $J: Y \times U \times \Gamma$ is a parameter-dependent objective functional,  $y\in Y$ and $u\in U$ denote the state and control functions respectively, $Y$ and $U$ are appropriate function spaces defined on $\Omega(\bxi)$. Both $y$ and $u$ depend on~$\bx$ and~$\bxi$, and the operator $\pF$ represents the governing equations. In this work, it is formulated as a system of parameter-dependent PDEs and consists of the partial differential equation operator $\pF_I$ and the boundary condition operator $\pF_B$. The set~$U_{ad}(\bxi)$ represents the admissible set, which is a parameter-dependent, bounded, closed, and convex subset of~$U$. It imposes an additional inequality constraint on $u(\bx,\bxi)$.

Since the optimal control system \eqref{OCP} is a constrained minimization problem, the necessary condition for the minimizer of \eqref{OCP} can be described by the following Karush-Kuhn-Tucker (KKT) system~\cite{de2015numerical, hinze2008optimization, troltzsch2010optimal}:
  \begin{equation}
    \left\{
    \begin{aligned}  
        &J_y\big(y^*(\bx, \bxi), u^*(\bx, \bxi); \bxi \big) -  \pF_y^*\big(y^*(\bx, \bxi), u^*(\bx, \bxi); \bxi \big)\, p^*(\bx, \bxi) = 0, \\
        &\pF\big(y^*(\bx, \bxi), u^*(\bx, \bxi); \bxi \big) = 0, \\   
        &(\mathrm{d}_uJ\big(y^*(\bx, \bxi), u^*(\bx, \bxi); \bxi \big), v(\bx, \bxi) - u^*(\bx, \bxi)) \geq 0, \quad \forall v(\bx, \bxi) \in U_{ad}(\bxi), \label{KKT}
    \end{aligned}
    \right.
  \end{equation}
where $y^*(\bx, \bxi), u^*(\bx, \bxi)$ represents the minimizer of \eqref{OCP}, and $p^*(\bx, \bxi)$ is the adjoint function (Lagrange multiplier) associated with the system. $\pF_y^*$ denotes the adjoint operator of $\pF_y$. Since $y(\bx,\bxi)$ is detremined by $u(\bx,\bxi)$ through the governing equation $\pF$, the derivative of the objective functional $J$ with respect to the control variable $u(\bx,\bxi)$ in the KKT system~\eqref{KKT} can be written as follows:
\begin{equation}\label{duj}
    \mathrm{d}_uJ\big(y^*(\bx, \bxi), u^*(\bx, \bxi); \bxi \big) = J_u\big(y^*(\bx, \bxi), u^*(\bx, \bxi); \bxi\big) - \pF_u^*\big(y^*(\bx, \bxi), u^*(\bx, \bxi); \bxi\big)p^*(\bx, \bxi). 
\end{equation}

The solution of the above OCP~\eqref{OCP} satisfies the KKT system~\eqref{KKT}. Therefore, instead of solving the primal problem~\eqref{OCP} directly, we seek solutions by solving the KKT system~\eqref{KKT}. However, directly solving the KKT system remains challenging due to the presence of parameter-dependent PDE operators and the additional control constraint. In this work, we propose a deep neural network surrogate modeling approach to address this issue. Our method employs separate neural networks to approximate the state function $y$, control function $u$, and adjoint function $p$, while incorporating an adaptive sampling strategy to mitigate the inefficiencies of random sampling—particularly for problems with low regularity or high-dimensional parameter spaces.

\subsection{The adjoint-oriented neural network method}
We begin by focusing on constructing surrogate models to solve the KKT system mentioned above. Let $\hat{y}(\bx, \bxi;\theta_y)$, $\hat{p}(\bx, \bxi;\theta_p)$, and $\hat{u}(\bx, \bxi;\theta_u)$ denote three distinct deep neural networks parameterized by $\theta_y$, $\theta_p$, and $\theta_u$, respectively. To address the parameter-dependent partial differential equations, we augment the input of each neural network by incorporating the model parameters $\bxi$ alongside the spatial variable $\bx$. This strategy aligns with approaches used for solving control problems~\cite{sun2020surrogate} and parameterized forward problems~\cite{khodayi2020varnet}. The input to these deep neural networks consist of both the spatial variable $\bx$ and parameters $\bxi$, which are concatenated as:
\begin{equation*}
    (\bx,\bxi) = [x^{(1)}, x^{(2)},\cdots, x^{(n)}, \xi^{(1)},\xi^{(2)}, \cdots, \xi^{(d)}],
\end{equation*}
where $\bx = [x^{(1)}, x^{(2)}, \cdots, x^{(n)}]$ represents the spatial coordinates, and $\bxi = [\xi^{(1)}, \xi^{(2)}, \cdots, \xi^{(d)}]$ denotes the parameter vector. A classical sampling strategy involves separately selecting collocation points $\{\bx_i\}$ from $\Omega$ and $\{\bxi_j\}$  from  $\Gamma$, respectively, and then forming combined sets $(\bx_i, \bxi_j) \in \Omega \times \Gamma$. However, in this work, due to the parameter-dependent nature of the $\Omega(\bxi)$ in geometry-parametric problems, the collocation points are sampled in pairs from the joint spatioparametric domain $\Omega_{\Gamma}$, where
\begin{equation*}
    \Omega_{\Gamma} = \{(\bx,\bxi) \,|\, \bx \in \Omega(\bxi), \ \bxi\in \Gamma\}.
\end{equation*}

AONN method employs a penalty-free technique \cite{lagaris1998artificial,sheng2021pfnn} to rigorously enforce the boundary conditions. Taking the control function $\hat{u}(\bx, \bxi;\theta_u)$ as an example, the core of the penalty-free trick is to introduce two additional neural networks to approximate $\hat{u}$. One neural network $\hat{u}_B$ is a priori satisfies the essential boundary conditions, and the other neural network $\hat{u}_I$ approximates the solution in the interior of the computational domain. This method eliminates training difficulties related to boundary condition enforcement, thereby improving accuracy and robustness in handling computational domains with complex geometries. The decomposition is constructed by 
\begin{equation*}
    \hat{u}(\bx, \bxi;\theta_u) = \hat{u}_B(\bx, \bxi;\theta_{u_B}) + l(\bx,\bxi) \hat{u}_I(\bx, \bxi;\theta_{u_I}),
\end{equation*}
where $\theta_u = \{\theta_{u_B}, \theta_{u_I}\}$, $\theta_{u_B}$ and $\theta_{u_I}$ refer to the parameters of two neural networks, respectively. $l(\bx,\bxi)$ is a necessary length factor function that satisfies the following properties:
\begin{equation*}
        \left\{\begin{aligned}
            \quad l(\bx,\bxi) = 0 \qquad  &\text{on} \;\,\partial \Omega(\xi), \\
            \quad l(\bx,\bxi) > 0 \qquad  &\text{in} \;\,\Omega(\xi).\\
         \end{aligned}\right.
\end{equation*}

Three deep neural networks $\hat{y}(\bx, \bxi;\theta_y)$, $\hat{p}(\bx, \bxi;\theta_p)$, and $\hat{u}(\bx, \bxi;\theta_u)$ are trained to approximate $y^*(\bx, \bxi), p^*(\bx, \bxi), u^*(\bx, \bxi)$ respectively, with the corresponding loss functions defined as follows:
\begin{align}
    J_s(\theta_y, \theta_u) &= \int_{\Omega_{\Gamma}}  r_s^2 \big(\hat{y}(\bx, \bxi;\theta_y), \hat{u}(\bx, \bxi;\theta_u); \bxi \big) \, \mathrm{d}\bx \mathrm{d}\bxi,  \label{eq:resfun1}\\
    J_a(\theta_y,\theta_u,\theta_p) &= \int_{\Omega_{\Gamma}}  r_a^2 \big(\hat{y}(\bx, \bxi;\theta_y), \hat{u}(\bx, \bxi;\theta_u),\hat{p}(\bx, \bxi;\theta_p); \bxi \big) \, \mathrm{d}\bx \mathrm{d}\bxi, \\
    J_u(\theta_u, u_{\text{step}}) &= \int_{\Omega_{\Gamma}}  r_u^2 \big(\hat{u}(\bx, \bxi;\theta_u), u_{\text{step}}(\bx, \bxi); \bxi \big) \, \mathrm{d}\bx \mathrm{d}\bxi, \label{eq:resfun3}
\end{align}
where the residual functionals are given by 
\begin{align}
    r_s \big(y(\bx, \bxi), u(\bx, \bxi); \bxi \big) &= \pF \big(y(\bx, \bxi), u(\bx,\bxi); \bxi),\label{eq:res_1}\\
    r_a \big(y(\bx, \bxi), u(\bx, \bxi),p(\bx, \bxi); \bxi \big) &= J_y \big(y(\bx, \bxi), u(\bx, \bxi); \bxi \big)- \pF_y^* \big(y(\bx, \bxi), u(\bx, \bxi); \bxi\big)  p(\bx,\bxi) ,\label{eq:res_2}\\
    r_u\big(u(\bx,\bxi),u_{\text{step}}(\bx,\bxi);\bxi\big)&=u(\bx,\bxi)-u_{\text{step}}(\bx,\bxi). \label{eq:res_3}
\end{align}

The term $u_{\text{step}}(\bx,\bxi)$ appearing in~\eqref{eq:resfun3} and~\eqref{eq:res_3} denotes an intermediate variable introduced to handle the additional control constraint. It plays a role in the iterative procedure for updating the control function in the third variational inequality of the KKT system~\eqref{KKT}, which is defined as follows:
\begin{equation}
    u_{\text{step}}(\bx, \bxi) = \mathbf{P}_{U_{\text{ad}}(\bxi)}\bigg(\hat{u}(\bx, \bxi,\theta_u)-c\mathrm{d}_uJ(\hat{y}(\bx, \bxi,\theta_y),\hat{u}(\bx,\bxi,\theta_u); \bxi)\bigg), 
\end{equation}
where $c>0$ denotes a prescribed step size, and $\mathbf{P}_{U_{\text{ad}}(\bxi)}$
denotes the projection operator onto the admissible control set $U_{\text{ad}}(\bxi)$. To illustrate the projection operator, consider a box constraint imposed on the control $u(\bx,\bxi)$: 
\begin{equation*}
    U_{ad}(\bxi) = \{ u \in U : u_a(\bx, \bxi) \leq u(\bx, \bxi)  \leq u_b(\bx, \bxi) \; \forall \bx \in \Omega(\bxi)\},
\end{equation*}
where $u_a$ and $u_b$ represent the lower and upper bound functions, respectively. Both $u_a$ and $u_b$ depend on parameter $\bxi$. Let $[u^{(1)},\cdots,u^{(N)}]^T$ denote the values of the control function at $N$ collocation points $\{x^{(i)},\xi^{(i)}\}_{i=1}^N$ within $\Omega_{\Gamma}$. The projection operator $\mathbf{P}_{U_{\text{ad}}(\bxi)}$ is given by~\cite{hinze2008optimization, troltzsch2010optimal}:
\begin{equation}\notag
        \mathbf{P}_{U_{\text{ad}}(\bxi)}(u^{(i)}) = 
        \begin{cases}
            u_a(\bx^{(i)},\bxi^{(i)}) \quad & \text{if}  \ u^{(i)} < u_a(\bx^{(i)}, \bxi^{(i)}), \\
        u^{(i)} \quad & \text{if} \ u_a(\bx^{(i)},\bxi^{(i)}) \leq u^{(i)} \leq u_b(\bx^{(i)},\bxi^{(i)}), \\
        u_b(\bx^{(i)},\bxi^{(i)}) \quad & \text{if}  \ u^{(i)} > u_b(\bx^{(i)}, \bxi^{(i)}), 
        \end{cases}
\end{equation}
for $i = 1,\cdots, N$. Note that from~\eqref{duj}, the total derivative $\mathrm{d}_uJ$ depends on the adjoint function $p^*(\bx,\bxi)$, and therefore, the loss function $J_u(\theta_u, u_{\text{step}})$ in \eqref{eq:resfun3} is influenced by $\theta_u$ , $\theta_y$ and $\theta_p$ together. By minimizing $J_u(\theta_u, u_{\text{step}})$ with respect to $\theta_u$, we can obtain an approximation of the optimal control $\hat{u}(\bx, \bxi;\theta_u)$.

In AONN method, neural networks are trained in an iterative manner. In each iteration, the network parameters are updated sequentially: $\hat{y}(\bx, \bxi;\theta_y)$ is trained first, followed by an update of $\hat{p}(\bx, \bxi;\theta_p)$, and finally refining $\hat{u}(\bx, \bxi;\theta_u)$. Specifically, we start with initial neural networks $\hat{y}(\bx, \bxi;\theta_y^{(0)})$, $\hat{p}(\bx, \bxi;\theta_p^{(0)})$, and $\hat{u}(\bx, \bxi;\theta_u^{(0)})$. At the $(i-1)$-th itertaion, neural networks are denoted by $\hat{y}(\bx, \bxi;\theta_y^{(i-1)})$, $\hat{p}(\bx, \bxi;\theta_p^{(i-1)})$, and $\hat{u}(\bx, \bxi;\theta_u^{(i-1)})$. Then, at the $i$-th iteration, the networks are updated as follows:
\begin{equation}\notag
\begin{split}
\text { training } \hat{y}: &\ \theta_y^{(i)}=\arg \min _{\theta_y} J_s\left(\theta_y, \theta_u^{(i-1)}\right), \\
\text { updating } \hat{p}: &\ \theta_p^{(i)}=\arg \min _{\theta_p} J_a\left(\theta_y^{(i)}, \theta_u^{(i-1)}, \theta_p\right), \\
\text { refining } \hat{u}: &\ \theta_u^{(i)}=\arg \min _{\theta_u} J_u\left(\theta_u, u_{\text {step }}^{(i-1)}\right).
\end{split}
\end{equation}
where
\begin{equation} \label{projection step}
    u_{\text{step}}^{(i-1)}(\bx, \bxi) = \mathbf{P}_{U_{\text{ad}}(\bxi)}\bigg(\hat{u}(\bx, \bxi;\theta_u^{(i-1)})-c^{(i)}\mathrm{d}_uJ(\hat{y}(\bx, \bxi;\theta_y^{(i)}),\hat{u}(\bx,\bxi;\theta_u^{(i-1)}); \bxi)\bigg), 
\end{equation}
and
\begin{equation}
\begin{split} \label{duj expression}
    \mathrm{d}_uJ\big(\hat{y}(\bx, \bxi;\theta_y^{(i)}), \hat{u}(\bx, \bxi;\theta_u^{(i-1)}); \bxi \big) & = J_u\big(\hat{y}(\bx, \bxi;\theta_y^{(i)}), \hat{u}(\bx, \bxi;\theta_u^{(i-1)}); \bxi\big) \\
    & - \pF_u^*\big(\hat{y}(\bx, \bxi;\theta_y^{(i)}), \hat{u}(\bx, \bxi;\theta_u^{(i-1)}); \bxi\big)\hat{p}(\bx, \bxi;\theta_p^{(i)}). 
\end{split}
\end{equation}

We summarize the process of constructing surrogate models to solve the optimal control problem in Algorithm~\ref{alg:AONN}. Note that before optimizing the loss functions $J_s(\theta_y)$, $J_a(\theta_p)$ and $J_u(\theta_u)$ with respect to $\theta_y$, $\theta_u$ and $\theta_p$, the integral defined in~\eqref{eq:resfun1}--\eqref{eq:resfun3} must undergo numerical discretization, leading to the following empirical loss: 
\begin{align}
        L_s(\theta_y, \theta_u) &= \bigg(\frac{1}{N}\sum_{i=1}^{N}\left|r_s \big(\hat{y}(\bx^{(i)}, \bxi^{(i)}; \theta_y), \hat{u}(\bx^{(i)}, \bxi^{(i)}; \theta_u); \bxi^{(i)} \big)\right|^2\bigg)^{1/2},\label{eq:res_num1} \\    
        L_a(\theta_y, \theta_u, \theta_p) &= \bigg(\frac{1}{N}\sum_{i=1}^{N}\left|r_a \big(\hat{y}(\bx^{(i)}, \bxi^{(i)}; \theta_y), \hat{u}(\bx^{(i)}, \bxi^{(i)}; \theta_u),\hat{p}(\bx^{(i)},\bxi^{(i)};\theta_p); \bxi^{(i)} \big)  \right|^2\bigg)^{1/2},\\
        L_u(\theta_u, u_{\text{step}}) &= \bigg(\frac{1}{N}\sum_{i=1}^{N}\left|r_u \big(\hat{u}(\bx^{(i)}, \bxi^{(i)}; \theta_u),u_{\text{step}}(\bx^{(i)}, \bxi^{(i)}); \bxi^{(i)}\big)\right|^2\bigg)^{1/2}, \label{eq:res_num3}
    \end{align}
where $\{(\bx^{(i)}, \bxi^{(i)})\}_{i=1}^N$ denotes the set of collocation points, which are typically sampled randomly from either a Gaussian or a uniform distribution.

The parameters of the neural networks can be updated efficiently using the automatic differentiation capabilities provided by deep learning libraries such as PyTorch \cite{paszke2017automatic} or TensorFlow \cite{abadi2016tensorflow}. For more details, readers are referred to~\cite{yin2024aonn}. 

\begin{breakablealgorithm}
\label{AONN_algorithm}
     \caption{AONN for solving parametric OCP($\bx, \bxi$)}\label{alg:AONN}
     \renewcommand{\algorithmicrequire}{\textbf{Input:}}
     \renewcommand{\algorithmicensure}{\textbf{Output:}}
     \begin{algorithmic}[1]
        \REQUIRE Initial network parameters $\theta_y^{(0)}, \theta_u^{(0)}, \theta_p^{(0)}$, decay factor $\gamma \in (0,1]$, initial step size $c^{(0)}$, initial number of epochs $n^{(0)}$, positive integer $n_{\text{aug}}$, maximum epoch number $N_{\text{ep}}$, batch size $m$, training set $\mathrm{S}_{\Omega}=\{\bx^{(i)}, \bxi^{(i)}\}_{i=1}^{n_r}$.
            \STATE // Train surrogate models
            \FOR{$i=1: N_{\text{ep}}$}
                \FOR{$j$ steps}
                    \STATE Sample $m$ samples from $\mathrm{S}_{\Omega}$.
                    \STATE $\theta_y^{(i)} \longleftarrow \text{arg} \min_{\theta_y} L_s(\theta_y, \theta_u^{(i-1)})$: Train network $\hat{y}(\bx, \bxi; \theta_y)$ with initialization $\theta_y^{(i-1)}$ for $n^{(i)}$ epochs.
                    \STATE $\theta_p^{(i)} \longleftarrow \text{arg} \min_{\theta_p} L_a(\theta_y^{(i)}, \theta_u^{(i-1)}, \theta_p)$: Train network $\hat{p}(\bx, \bxi; \theta_p)$ with initialization $\theta_p^{(i-1)}$ for $n^{(i)}$ epochs.
                    \STATE Compute $u_{\text{step}}^{(i-1)}(\bx, \bxi)$ by \eqref{projection step} and \eqref{duj expression} with $\theta_y^{(i)}$, $\theta_p^{(i)}$ and $\theta_u^{(i-1)}$.
                    \STATE $\theta_u^{(i)} \longleftarrow \text{arg} \min_{\theta_u} L_u(\theta_u, u_{\text{step}}^{(i-1)} )$: Train network $\hat{u}(\bx, \bxi; \theta_u)$ with initialization $\theta_u^{(i-1)}$ for $n^{(i)}$ epochs.
                    \STATE $c^{(i+1)} = \gamma c^{(i)}$.
                    \STATE $n^{(i+1)} = n^{(i)} + n_{\text{aug}}$.
                    \STATE Let $\theta_y^{*} = \theta_y^{(i)}$ and $\theta_p^{*} = \theta_p^{(i)}$ and $\theta_u^{*} = \theta_u^{(i)}$.
                \ENDFOR
            \ENDFOR
        \ENSURE $\hat{y}(\bx, \bxi; \theta_y^{*}), \hat{u}(\bx, \bxi; \theta_u^{*})$ 
    \end{algorithmic}
\end{breakablealgorithm}

\section{Deep adaptive sampling for surrogate modeling} \label{section_sampling}
In recent years, the concept of adaptive sampling has attracted increasing attention, leading to the emergence of numerous adaptive sampling techniques tailored to diverse scenarios. In the context of the problem addressed in this paper, solving parametric OCPs constrained by PDEs, adaptive methods can similarly be employed to guide the acquisition of collocation points required for surrogate modeling.

Compared to uniform sampling, collocation points obtained via adaptive strategies make a more substantial contribution to the training process. Consequently, these methods not only enhance the training speed and accuracy of the surrogate model, but also effectively address complex challenges that simple uniform sampling fails to handle, such as high-dimensional problems or those involving singularities. In this section, we provide a brief introduction to the deep adaptive sampling method for surrogate modeling without labeled data ($\text{DAS}^2$)~\cite{wang2024deep} method, which is particularly well suited for the construction of surrogate models.

The key idea of $\text{DAS}^2$ is to generate random collocation points based on the residual-induced distribution, with the aim of reducing the statistical error introduced by discretizing the loss functional using Monte Carlo approximation. Specifically, in this work, the residual $r(\bx,\bxi; \theta)$ is given by:

\begin{equation}
 r^2(\bx,\bxi; \theta) =  r_s^2 \big(\hat{y}(\bx, \bxi;\theta_y), \hat{u}(\bx, \bxi;\theta_u); \bxi \big) + r_a^2 \big(\hat{y}(\bx, \bxi;\theta_y), \hat{u}(\bx, \bxi;\theta_u),\hat{p}(\bx, \bxi;\theta_p); \bxi \big) , \label{residual_eqn}
\end{equation}%
where $\theta = (\theta_y, \theta_p, \theta_u)$ represents the parameters of three neural networks, and $r_s$ and $r_a$ are defined in \eqref{eq:res_1} and \eqref{eq:res_2} respectively.
In $\text{DAS}^2$, we adaptively sample from a joint probability density function (PDF), $\hat{r}(\bx, \bxi)$, induced by the residual $r(\bx,\bxi; \theta)$ for a given $\theta$:
\begin{equation}\label{eq:pdf_res}
    \hat{r}(\bx, \bxi) \propto r^2(\bx,\bxi; \theta) \cdot h(\bx, \bxi),
\end{equation}
where $h(\bx, \bxi)$ is a cutoff function defined on a compact support $B \supset \Omega(\bxi) \times \Gamma$. Here, $B$ is chosen to be slightly larger than $\Omega(\bxi) \times \Gamma$. In this work, we adopt the cutoff function defined in \cite{tang2023pinns}, where $h(\bx, \bxi) = 1$ if $(\bx, \bxi) \in \Omega(\bxi) \times \Gamma$ and linearly decays to $0$ towards $\partial B$. 

In $\text{DAS}^2$, the deep generative model, called KRnet \cite{tang2020deep,tang2022adaptive,tang2023pinns}, is employed to approximate the residual-induced distribution and generate the corresponding random collocation points. The PDF induced by KRnet is given by:
\begin{equation}
    p_{\text{KRnet}}(\bx, \bxi; \theta_f) = p_{\boldsymbol{Z}}\big(f_{\text{KRnet}}(\bx, \bxi; \theta_f)\big) \left|\text{det}\nabla_{\bx, \bxi}f_{\text{KRnet}} \right|, \label{KRnet}
\end{equation}
where $f_{\text{KRnet}}(\bx, \bxi; \theta_f)$ is a specially designed invertible mapping parameterized with neural networks, and $p_{\boldsymbol{Z}}$ denotes the prior distribution of the random vector $\boldsymbol{z}$.

The structure of KRnet is the core of $\text{DAS}^2$. Here, we briefly outline its workflow:
\begin{equation*}
    \boldsymbol{z} = f_{\text{KRnet}}(\bx, \bxi;\theta_f) = L_N \circ f_{K-1}^{\text{outer}} \circ \cdots \circ f_{1}^{\text{outer}}(\bx, \bxi; \theta_f),
\end{equation*}
where $f_{k}^{\text{outer}}$ denotes the $k$-th iteration of the outer loop, defined as:
\begin{equation*}
    f_{k}^{\text{outer}} = L_S \circ f_{[k,L]}^{\text{inner}} \circ \cdots \circ f_{[k,1]}^{\text{inner}} \circ L_R,
\end{equation*}
where $f_{[k,j]}^{\text{inner}}$ refers to a combination of $L$ affine coupling layers, each followed by a scale layer and bias layer, as proposed in \cite{kingma2018glow, dinh2016density}. The layers $L_N$, $L_S$, and $L_R$ correspond to the nonlinear, squeezing, and rotation layers, respectively. The prior distribution $p_Z$ for the random vector $\boldsymbol{z}$ is usually chosen as the standard normal distribution.

To generate random samples that are consistent with $\hat{r}(\bx,\bxi)$, we need to find optimal parameters $\theta_f^*$ such that the PDF model $p_{\text{KRnet}}(\bx, \bxi; \theta_f^*)$ can approximate $\hat{r}(\bx,\bxi)$ most. This task can be reformulated as minimizing the following cross-entropy loss $H(\hat{r}, p_{\text{KRnet}})$ \cite{de2005tutorial}.
\begin{equation}
    H(\hat{r}, p_{\text{KRnet}}) \approx -\frac{1}{M}\sum_{i=1}^{M}\frac{\hat{r}(\bx^{(i)},\bxi^{(i)})\log p_{\text{KRnet}}(\bx^{(i)},\bxi^{(i)};\theta_f)}{p_{\text{KRnet}}(\bx^{(i)},\bxi^{(i)};\hat{\theta}_f)}, \label{importance sampling}
\end{equation}
where $M$ is the number of collocation points and $\hat{\theta}_f$ is the parameters of a known PDF model $p_{\text{KRnet}}(\bx^{(i)},\bxi^{(i)};\hat{\theta}_f)$. We will specify how to choose $p_{\text{KRnet}}(\bx^{(i)},\bxi^{(i)};\hat{\theta}_f)$ in Section~\ref{sec:adaptivewithdas}. Once the optimal parameters $\theta_f^*$ are determined by minimizing \eqref{importance sampling}, samples from $p_{\text{KRnet}}(\bx, \bxi; \theta_f^*)$ can be obtained through the inverse mapping of KRnet:
$$(\bx, \bxi) = f_{\text{KRnet}}^{-1}(\boldsymbol{z}; \theta_f^*),$$
by sampling $\boldsymbol{z}$. 
We can then use these obtained samples to train the surrogate model mentioned above. Note that only those samples within $\Omega(\bxi)\times\Gamma$ are retained. We summarize the process of $\text{DAS}^2$ for training a surrogate model $\hat{s}(\bx,\bxi;\theta)$ with loss function $L(\bx,\bxi;\theta)$ in Algorithm~\ref{alg:DAS}. For more details, readers are referred to~\cite{tang2020deep,tang2022adaptive,tang2023pinns}.

\begin{breakablealgorithm}
\label{DAS_algorithm}
     \caption{$\text{DAS}^2$ for training a surrogate model}\label{alg:DAS}
     \renewcommand{\algorithmicrequire}{\textbf{Input:}}
     \renewcommand{\algorithmicensure}{\textbf{Output:}}
     \begin{algorithmic}[1]
        \REQUIRE Loss function $L(\theta)$, initial $p_{\mathrm{KRnet}}(\bx, \bxi ; \theta_f^{(0)})$, initial network parameters $\theta^{(0)}$, maximum epoch number $N_{\text{ep}}$, batch size $m$, initial training set $\mathrm{S}_{\Omega, 0}=\{\bx^{(i)}, \bxi^{(i)}\}_{i=1}^{n_r}$.
        \FOR{$k=0: N_{\text {adaptive}}-1$} 
            \STATE // Train surrogate models
            \FOR{$i=1: N_{\text{ep}}$}
                \FOR{$j$ steps}
                    \STATE Sample $m$ samples from $\mathrm{S}_{\Omega, k}$.
                    \STATE 
                    $\theta^{(i)} \longleftarrow \text{arg} \min_{\theta} L(\theta)$: Train network $\hat{s}(\bx, \bxi; \theta)$ with initialization $\theta^{(i-1)}$.
                    \STATE Update  $\hat{r}(\bx,\bxi)$ according to \eqref{eq:pdf_res}.
                    \STATE Let $\theta^{*} = \theta^{(i)}$.
                \ENDFOR
            \ENDFOR
            \STATE Let $\theta^{(0)} = \theta^{*}$.
            \STATE // Update KRnet
            \FOR{$i=1: N_{\text{ep}}$}
                \FOR{$j$ steps}
                    \STATE Sample $m$ samples from $p_{\text {KRnet }}\left(\boldsymbol{x}, \bxi ; \theta_f^{*,(k-1)}\right)$.
                    \STATE Update $p_{\text {KRnet }}\left(\boldsymbol{x}, \bxi ; \theta_f\right)$ by descending the stochastic gradient of $H\left(\hat{r}, p_{\text {KRnet }}\right)$. 
                \ENDFOR
            \ENDFOR
            \STATE // Refine training set
            \STATE Generate $\mathrm{S}_{\Omega, k+1}^g \subset \Omega_{\Gamma}$ with size $n_r$ through $p_{\text {KRnet}}\left(\boldsymbol{x}, \bxi ; \theta_f^{*,(k+1)}\right)$.
            \STATE $\mathrm{S}_{\Omega, k+1}=\mathrm{S}_{\Omega, k} \cup \mathrm{S}_{\Omega, k+1}^g$.
        \ENDFOR
        \ENSURE Surrogate model $s(\bx, \bxi; \theta^{*})$
    \end{algorithmic}
\end{breakablealgorithm}

\section{Adaptive AONN for parametric optimal control problems}\label{section_method}

In this section, we propose the adaptive AONN for solving parametric optimal control problems, which integrates AONN with adaptive sampling techniques. In particular, $\text{DAS}^2$ is embedded into the AONN training process to enhance computational robustness and generalization performance.

\subsection{Adaptive AONN} 
In the previous section, we described how the construction of the surrogate model allows us to circumvent the direct solution of the original complex problem and instead solve the OCP through its associated KKT system. In addition to building the surrogate model, the choice of an appropriate sampling strategy is also essential. It should be noted that training surrogate models using collocation points sampled from a uniform distribution often yields poor performance, and it has been observed that the AONN approach exhibits notable inaccuracies when solving low-regularity problems~\cite{tang2022adaptive, tang2023pinns, wu2023comprehensive}. Due to the large variance in low-regularity solutions, most collocation points sampled from a uniform distribution contribute little to the approximation. Since both the spatial and parametric domains need to be sampled, the additional dimensions introduced by the parametric domain~$\Gamma$ further complicate the sampling process~\cite{wright2022high}. 

To overcome these challenges, we introduce the adaptive sampling strategy into the framework to replace random sampling and thereby construct an adaptive AONN method. The idea behind this approach is straightforward: we alternately update the surrogate models and their training sets in an adaptive manner. At each stage, we first train the surrogate models by minimizing the empirical loss functions defined in \eqref{eq:res_num1}–\eqref{eq:res_num3}, using collocation points sampled according to the residual-induced distribution defined in \eqref{eq:pdf_res}. We then apply the newly optimized surrogate models to update the residual-induced distribution. Adaptive sampling significantly improves sampling efficiency, and by selecting a suitable adaptive strategy, we ensure that newly sampled collocation points in each iteration provide maximal benefit to the next stage of training.

After incorporating the adaptive sampling procedure into the DAL framework, the intended workflow can be formulated in Algorithm~\ref{alg:aAONN}.

\begin{breakablealgorithm}
\label{aAONN_algorithm}
     \caption{Adaptive AONN for solving parametric OCP($\bx, \bxi$)}\label{alg:aAONN}
     \renewcommand{\algorithmicrequire}{\textbf{Input:}}
     \renewcommand{\algorithmicensure}{\textbf{Output:}}
     \begin{algorithmic}[1]
        \REQUIRE Initial network parameters $\theta_y^{(0)}, \theta_u^{(0)}, \theta_p^{(0)}$, decay factor $\gamma \in (0,1]$, initial step size $c^{(0)}$, initial number of epochs $n^{(0)}$, positive integer $n_{\text{aug}}$, maximum epoch number $N_{\text{ep}}$, batch size $m$, initial training set $\mathrm{S}_{\Omega, 0}=\{\bx^{(i)}, \bxi^{(i)}\}_{i=1}^{n_r}$.
        \FOR{$k=0: N_{\text {adaptive}}-1$} 
            \STATE // Train surrogate models
            \FOR{$i=1: N_{\text{ep}}$}
                \FOR{$j$ steps}
                    \STATE Sample $m$ samples from $\mathrm{S}_{\Omega, k}$.
                    \STATE $\theta_y^{(i)} \longleftarrow \text{arg} \min_{\theta_y} L_s(\theta_y, \theta_u^{(i-1)})$: Train network $\hat{y}(\bx, \bxi; \theta_y)$ with initialization $\theta_y^{(i-1)}$ for $n^{(i)}$ epochs.
                    \STATE $\theta_p^{(i)} \longleftarrow \text{arg} \min_{\theta_p} L_a(\theta_y^{(i)}, \theta_u^{(i-1)}, \theta_p)$: Train network $\hat{p}(\bx, \bxi; \theta_p)$ with initialization $\theta_p^{(i-1)}$ for $n^{(i)}$ epochs.
                    \STATE Compute $u_{\text{step}}^{(i-1)}(\bx, \bxi)$ by \eqref{projection step} \eqref{duj expression} with $\theta_y^{(i)}$, $\theta_p^{(i)}$ and $\theta_u^{(i-1)}$.
                    \STATE $\theta_u^{(i)} \longleftarrow \text{arg} \min_{\theta_u} L_u(\theta_u, u_{\text{step}}^{(i-1)} )$: Train network $\hat{u}(\bx, \bxi; \theta_u)$ with initialization $\theta_u^{(i-1)}$ for $n^{(i)}$ epochs.
                    \STATE $c^{(i+1)} = \gamma c^{(i)}$.
                    \STATE $n^{(i+1)} = n^{(i)} + n_{\text{aug}}$.
                    \STATE Let $\theta_y^{*} = \theta_y^{(i)}$ and $\theta_p^{*} = \theta_p^{(i)}$ and $\theta_u^{*} = \theta_u^{(i)}$.
                \ENDFOR
            \ENDFOR
            \STATE Let $\theta_y^{(0)} = \theta_y^{*}$ and $\theta_p^{(0)} = \theta_p^{*}$ and $\theta_u^{(0)} = \theta_u^{*}$.
            \STATE // Sampling collocation points
            \STATE Update the residual-induced distribution \eqref{eq:pdf_res}, and apply adaptive sampling method to obtain $\mathrm{S}_{\Omega, k+1}$.
            \STATE $\mathrm{S}_{\Omega, k+1}=\mathrm{S}_{\Omega, k} \cup \mathrm{S}_{\Omega, k+1}^g$.
        \ENDFOR
        \ENSURE $\hat{y}(\bx, \bxi; \theta_y^{*}), \hat{u}(\bx, \bxi; \theta_u^{*})$ 
    \end{algorithmic}
\end{breakablealgorithm}

\subsection{Adaptive AONN with $\text{DAS}^2$}\label{sec:adaptivewithdas}

The practical realization of the adaptive AONN requires the selection of a concrete adaptive sampling strategy within the iterative training process. Since $\text{DAS}^2$ is capable of generating samples adaptively according to a residual-induced distribution, it is employed to construct the training set throughout the adaptive iterations. At each stage, we first train the surrogate models by minimizing the empirical loss functions defined in \eqref{eq:res_num1}–\eqref{eq:res_num3}, using collocation points sampled by KRnet obtained from the previous iteration. TThen,we apply the newly optimized surrogate models to update the parameters of KRnet by minimizing the cross-entropy loss function defined in \eqref{importance sampling}. Automatic differentiation tools provided by PyTorch \cite{paszke2017automatic} or TensorFlow \cite{abadi2016tensorflow} are utilized to compute the derivatives within those loss functions.

More specifically, the mechanism begins with the initialization of the deep generative model $p_{\mathrm{KRnet}}(\bx, \bxi ; \theta_f^{(0)})$, the network parameters $\theta_y^{(0)}, \theta_u^{(0)}$ and $\theta_p^{(0)}$, and the acquisition of the initial training set $\mathrm{S}_{\Omega, 0}$. We first draw $m$ samples from $\mathrm{S}_{\Omega, 0}$ required for the first adaptive training iteration. Subsequently, the parameters $\theta_y^{(0)}$ of the state network $\hat{y}(\bx, \bxi; \theta_y)$ are updated by minimizing the loss function $L_s(\theta_y, \theta_u^{(0)})$. Using the updated $\theta_y^{(1)}$, we then minimize the loss function $L_a(\theta_y^{(1)}, \theta_u^{(0)}, \theta_p)$ to update the parameters $\theta_p^{(0)}$ of the adjoint network $\hat{p}(\bx, \bxi; \theta_p)$. Following the updates of the state $\hat{y}$ and the adjoint $\hat{p}$, and prior to updating the control $\hat{u}$, the intermediate term $u_{\text{step}}^{(0)}$ is computed using $\theta_y^{(1)}$, $\theta_p^{(1)}$ and $\theta_u^{(0)}$ as defined in~\eqref{projection step} \eqref{duj expression}. Once $u_{\text{step}}^{(0)}$ is obtained, the parameters $\theta_u^{(0)}$ of the control network $\hat{u}(\bx, \bxi; \theta_u)$ are updated by minimizing the loss function $L_u(\theta_u, u_{\text{step}}^{(0)})$. Upon completion of the surrogate model training procedure, hyper-parameters including the step size $c$ and the number of epochs $n$ are updated. Here, the initial step size $c^{(0)}$ serves to ensure algorithm convergence, while the decay factor $\gamma$ is employed to enhance robustness. In each inner iteration, the number of epochs is incremented by $n_{\text{aug}}$; in contrast to a fixed epoch count, this progressively increasing number of epochs improves the accuracy of the algorithm. 

The aforementioned steps are repeated $N_{\text{ep}}$ times to conclude the training phase of the surrogate model within the first adaptive iteration, yielding the optimal network parameters $\theta_y^*$, $\theta_p^*$, and $\theta_u^*$. Subsequently, we utilize $p_{\text {KRnet}}\left(\boldsymbol{x}, \bxi ; \theta_f^{*,(0)}\right)$ = $p_{\mathrm{KRnet}}(\bx, \bxi ; \theta_f^{(0)})$ to draw $m$ samples required for the initial update of the sampling model. Substituting these samples and optimal network parameters into~\eqref{residual_eqn} yields the corresponding residual $r^2(\bx,\bxi; \theta^*)$. We select $p_{\mathrm{KRnet}}(\bx, \bxi ; \theta_f^{(0)})$ as the known PDF model on the right-hand side of~\eqref{importance sampling}. Finally, $\theta_f$ is updated by applying gradient descent to $H\left(\hat{r}, p_{\text {KRnet }}\right)$. Upon obtaining the updated deep generative model $p_{\text {KRnet}}\left(\boldsymbol{x}, \bxi ; \theta_f^{*,(1)}\right)$, we utilize it to generate the new collocation set $\mathrm{S}_{\Omega, 1}^g$ of size $n_r$. The training set $\mathrm{S}_{\Omega, 1}$ for the subsequent adaptive iteration is then constructed as the union of the initial training set $\mathrm{S}_{\Omega, 0}$ and the newly generated collocation set $\mathrm{S}_{\Omega, 1}^g$. 

This marks the completion of the first adaptive iteration. The subsequent iteration is initialized using the current optimal parameters $\theta_y^{*}$, $\theta_p^{*}$, $\theta_u^{*}$, and $\theta_f^{*,(1)}$. This adaptive process is repeated until a predetermined number of iterations $N_{\text {adaptive}}$ is reached. Upon termination, the algorithm outputs the optimal approximations for the state and control functions, denoted as $\hat{y}(\bx, \bxi; \theta_y^{*})$ and $\hat{u}(\bx, \bxi; \theta_u^{*})$, obtained from the final iteration.

\begin{figure}[h] 
    \centering
    \includegraphics[width=0.9\textwidth]{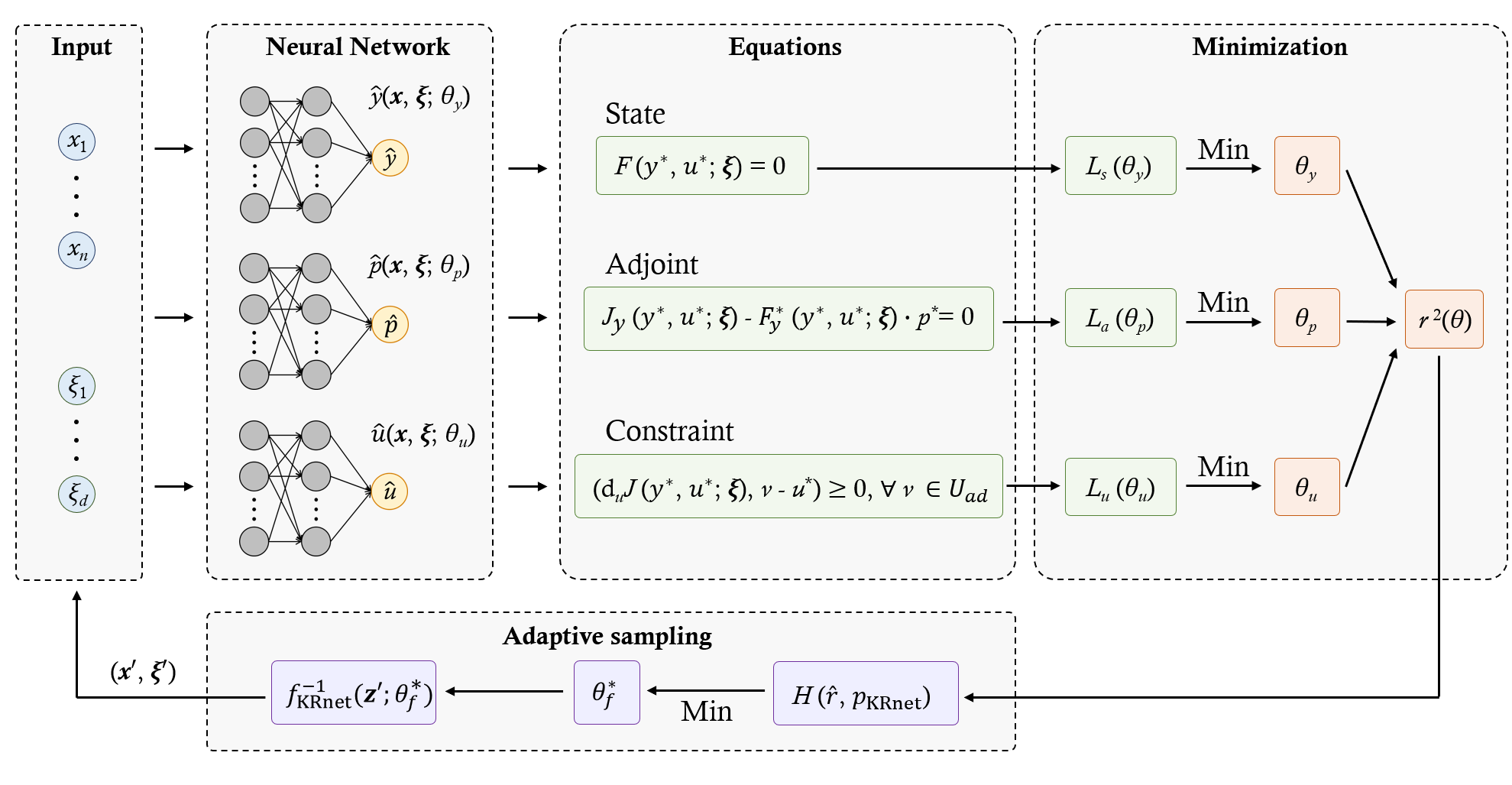}
    \caption{The workflow of adaptive AONN.}
    \label{workflow}
\end{figure}

In each stage, the solution obtained through AONN improves the ability of $\text{DAS}^2$ to approximate the target distribution, while samples generated by $\text{DAS}^2$ improve the training of the surrogate model in AONN, forming a mutually strengthening computational loop. Once the surrogate models are constructed and the deep generative model is trained, the corresponding optimal control function and state function can be efficiently evaluated for any set of parameters without the need to repeatedly solve the optimal control problem. The general procedure of the adaptive AONN with $\text{DAS}^2$ is summarized in Algorithm \ref{DAS-AONN_algorithm}, whereas the overall workflow is visualized in Fig \ref{workflow}.

\begin{breakablealgorithm}
\label{DAS-AONN_algorithm}
     \caption{adaptive AONN with $\text{DAS}^2$ for solving parametric OCP($\bx, \bxi$)}\label{alg:DAS-AONN}
     \renewcommand{\algorithmicrequire}{\textbf{Input:}}
     \renewcommand{\algorithmicensure}{\textbf{Output:}}
     \begin{algorithmic}[1]
        \REQUIRE Initial $p_{\mathrm{KRnet}}(\bx, \bxi ; \theta_f^{(0)})$, initial network parameters $\theta_y^{(0)}, \theta_u^{(0)}, \theta_p^{(0)}$, decay factor $\gamma \in (0,1]$, initial step size $c^{(0)}$, initial number of epochs $n^{(0)}$, positive integer $n_{\text{aug}}$, maximum epoch number $N_{\text{ep}}$, batch size $m$, initial training set $\mathrm{S}_{\Omega, 0}=\{\bx^{(i)}, \bxi^{(i)}\}_{i=1}^{n_r}$.
        \FOR{$k=0: N_{\text {adaptive}}-1$} 
            \STATE // Train surrogate models
            \FOR{$i=1: N_{\text{ep}}$}
                \FOR{$j$ steps}
                    \STATE Sample $m$ samples from $\mathrm{S}_{\Omega, k}$.
                    \STATE $\theta_y^{(i)} \longleftarrow \text{arg} \min_{\theta_y} L_s(\theta_y, \theta_u^{(i-1)})$: Train network $\hat{y}(\bx, \bxi; \theta_y)$ with initialization $\theta_y^{(i-1)}$ for $n^{(i)}$ epochs.
                    \STATE $\theta_p^{(i)} \longleftarrow \text{arg} \min_{\theta_p} L_a(\theta_y^{(i)}, \theta_u^{(i-1)}, \theta_p)$: Train network $\hat{p}(\bx, \bxi; \theta_p)$ with initialization $\theta_p^{(i-1)}$ for $n^{(i)}$ epochs.
                    \STATE Compute $u_{\text{step}}^{(i-1)}(\bx, \bxi)$ by \eqref{projection step} \eqref{duj expression} with $\theta_y^{(i)}$, $\theta_p^{(i)}$ and $\theta_u^{(i-1)}$.
                    \STATE $\theta_u^{(i)} \longleftarrow \text{arg} \min_{\theta_u} L_u(\theta_u, u_{\text{step}}^{(i-1)} )$: Train network $\hat{u}(\bx, \bxi; \theta_u)$ with initialization $\theta_u^{(i-1)}$ for $n^{(i)}$ epochs.
                    \STATE $c^{(i+1)} = \gamma c^{(i)}$.
                    \STATE $n^{(i+1)} = n^{(i)} + n_{\text{aug}}$.
                    \STATE Let $\theta_y^{*} = \theta_y^{(i)}$ and $\theta_p^{*} = \theta_p^{(i)}$ and $\theta_u^{*} = \theta_u^{(i)}$.
                \ENDFOR
            \ENDFOR
            \STATE Let $\theta_y^{(0)} = \theta_y^{*}$ and $\theta_p^{(0)} = \theta_p^{*}$ and $\theta_u^{(0)} = \theta_u^{*}$.
            \STATE // Update KRnet
            \FOR{$i=1: N_{\text{ep}}$}
                \FOR{$j$ steps}
                    \STATE Sample $m$ samples from $p_{\text {KRnet }}\left(\boldsymbol{x}, \bxi ; \theta_f^{*,(k)}\right)$.
                    \STATE Compute $r^2(\bx,\bxi; \theta^{*})$ by \eqref{residual_eqn} with $\theta_y^{*}$, $\theta_p^{*}$ and $\theta_u^{*}$.
                    \STATE $p_{\text {KRnet }}\left(\boldsymbol{x}, \bxi ; \hat \theta_f\right)$ = $p_{\text {KRnet }}\left(\boldsymbol{x}, \bxi ; \theta_f^{*,(k)}\right)$.
                    \STATE Update $p_{\text {KRnet }}\left(\boldsymbol{x}, \bxi ; \theta_f\right)$ by descending the stochastic gradient of $H\left(\hat{r}, p_{\text {KRnet }}\right)$. (see equation~\eqref{importance sampling})
                \ENDFOR
            \ENDFOR
            \STATE // Refine training set
            \STATE Generate $\mathrm{S}_{\Omega, k+1}^g \subset \Omega_{\Gamma}$ with size $n_r$ through $p_{\text {KRnet}}\left(\boldsymbol{x}, \bxi ; \theta_f^{*,(k+1)}\right)$.
            \STATE $\mathrm{S}_{\Omega, k+1}=\mathrm{S}_{\Omega, k} \cup \mathrm{S}_{\Omega, k+1}^g$.
        \ENDFOR
        \ENSURE $\hat{y}(\bx, \bxi; \theta_y^{*}), \hat{u}(\bx, \bxi; \theta_u^{*})$ 
    \end{algorithmic}
\end{breakablealgorithm}

Our proposed method, adaptive AONN with $\text{DAS}^2$, organically integrates the training framework of AONN with the sampling strategy of $\text{DAS}^2$. This integration not only overcomes the limitations of AONN—specifically its inability to solve problems exhibit singularities due to the lack of sampling technique optimization—but also addresses the precision issues inherent in $\text{DAS}^2$, which arise from utilizing fully connected neural networks for solving PDEs. Another key difference between $\text{DAS}^2$ and adaptive AONN with $\text{DAS}^2$ is that $\text{DAS}^2$ cannot handle additional constraints imposed on collocation points. It updates the sampling strategy solely on the basis of the residual distribution and thus cannot accommodate other potential restrictions on the sampling domain. In contrast, adaptive AONN with $\text{DAS}^2$ addresses this limitation by incorporating the DAL framework, thereby enabling adaptive sampling capable of handling inequality constraints. 

In the subsequent numerical experiments, we compare the performance of adaptive AONN with $\text{DAS}^2$ against AONN and $\text{DAS}^2$ in solving parametric OCPs constrained by PDEs subject to additional control constraints. Notably, the computational domains of these problems are characterized by singularities. The experimental results will demonstrate that our proposed method achieves a significant advantage in accuracy over existing state-of-the-art methods when addressing such complex problems involving singularities.

\section{Numerical study}\label{section_experiments}
In this section, we present three numerical experiments with singular features pronounced to validate the efficiency of the proposed method. To facilitate a comparative analysis, we evaluate the performance of adaptive AONN (with $\text{DAS}^2$), $\text{DAS}^2$ and AONN under identical problem configurations. High-fidelity solutions obtained via the finite element method are served as the reference solutions.

\subsection{Test 1: Optimal control for the Laplace equation with geometrical parametrization}
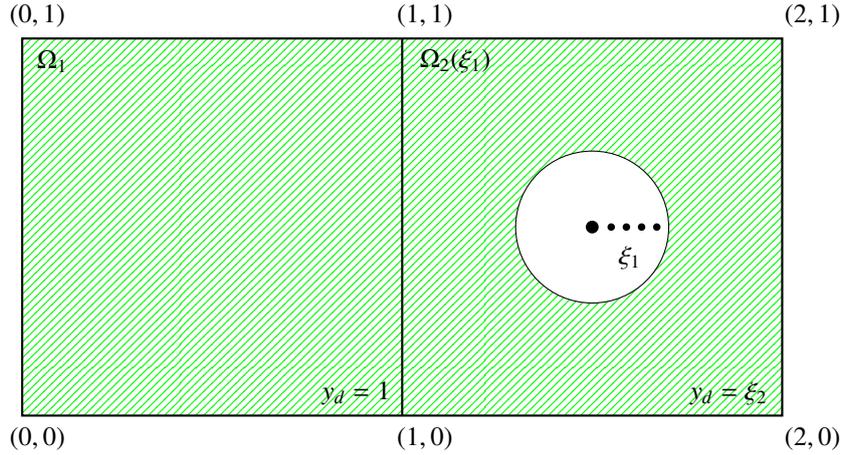
\begin{figure}[htbp]
    \centering
    \begin{tikzpicture}[thick]
        \draw[pattern=north east lines, pattern color=green, even odd rule] (0, 0) rectangle (10, 5);
        \draw[thick] (5, 0) -- (5, 5);
        \draw[thick] (7.5, 2.5) circle (1);
        \fill[white!20] (7.5,2.5) circle (1);
        \node at (0.2, -0.3) {$(0,0)$};
        \node at (0.2, 5.3) {$(0,1)$};
        \node at (5.3, -0.3) {$(1,0)$};
        \node at (10.4, -0.3) {$(2,0)$};
        \node at (5.3, 5.3) {$(1,1)$};
        \node at (10.4, 5.3) {$(2,1)$};    
        \node at (4.4, 0.3) {$y_d = 1$};
        \node at (9.3, 0.3) {$y_d = \xi_2$};
        \node at (0.4, 4.7) {$\Omega_1$};
        \node at (5.7, 4.7) {$\Omega_2(\xi_1)$};
        \filldraw[black] (7.5, 2.5) circle (2pt);
        \node at (8.0, 2.1) {$\xi_1$};
        \filldraw[black] (7.75, 2.5) circle (1pt);
        \filldraw[black] (7.95, 2.5) circle (1pt);
        \filldraw[black] (8.15, 2.5) circle (1pt);
        \filldraw[black] (8.35, 2.5) circle (1pt); 
    \end{tikzpicture}
    \caption{The parameter-dependent computational domain $\Omega$ in Test 1.}
    \label{fig:parameter_computation_domain}
\end{figure}

In this problem, we consider the following parametric OCP:
\begin{equation}\label{eq:test1}
\left\{\begin{aligned}
&\min \limits_{y(\bx, \bxi), u(\bx, \bxi)} J\big(y(\bx, \bxi), u(\bx, \bxi)\big)=\dfrac{1}{2}\left\|y(\bx,\bxi)-y_d(\bx,\bxi)\right\|_{2, \Omega}^2+\dfrac{\alpha}{2}\|u(\bx, \bxi)\|_{2, \Omega}^2, \\
&\text { subject to }\left\{\begin{aligned}
-\Delta y(\bx, \bxi)&=u(\bx, \bxi) && \text { in } \Omega, \\
y(\bx, \bxi)&=1 && \text { on } \partial \Omega,
\end{aligned}\right. \\ 
&\text { and } u_a \leq u(\bx, \bxi) \leq u_b \ \text { a.e. in } \Omega,
\end{aligned}\right.
\end{equation}
where $\alpha = 0.001$, $u_a = 0$, and $u_b = 10$. The vector $\bxi=\big(\xi_1, \xi_2\big)$ denotes the parameters, whose range is $\Gamma = [0.05, 0.45] \times [0.5,2.5]$. Both the spatial domain $\Omega$ and the desired state $y_d(\bx, \bxi)$ depend on $\bxi$. Specifically, the computational domain $\Omega$ is determined by the first parameter $\xi_1$ and is given by $\Omega=\big([0,2] \times[0,1]\big) \backslash \mathbb{B}\big((1.5,0.5), \xi_1\big)$, where~$\mathbb{B}\big((1.5,0.5),\xi_1\big)$ is a ball of radius $\xi_1$ centered at $(1.5,0.5)$. Figure~\ref{fig:parameter_computation_domain} illustrates the parameter-dependent computational domain. The desired state $y_d(\bx, \bxi)$ depends on the second parameter $\xi_2$ and is given by:
\begin{equation*}
    y_d(\bx, \bxi)= 
    \begin{cases}1 & \text { in } \Omega_1=[0,1] \times[0,1], \\
    \xi_2 & \text { in } \Omega_2=\big([1,2] \times[0,1]\big) \backslash \mathbb{B}\big((1.5,0.5), \xi_1\big).
    \end{cases}  
\end{equation*}

This test problem originates from a real-world application in the medical field~\cite{negri2013reduced, hesthaven2016certified}, which focuses on the treatment of cancer by local hyperthermia. In this application, the aim is to control heating sources to create distinct temperature zones: maintaining a higher temperature in the tumor area and a lower temperature in the surrounding non-lesion regions. The right side of Figure~\ref{fig:parameter_computation_domain} illustrates an organ that contains a tumor. 

In this context, our goal is to develop a surrogate model for optimal heat source control that works effectively for different organ shapes and temperature patterns resulting from various set of parameters. Variations in the geometric parameter $\xi_1$ lead to different computational domains, posing challenges for conventional mesh-based numerical methods that struggle to handle irregular boundary adaptation efficiently. However, this issue can be addressed through parametric space sampling using neural networks. The spatioparametric space considered in this example is defined as:
\begin{equation*}
    \begin{aligned}
    \Omega_{\Gamma} \coloneqq \{ (x_1, x_2, \xi_1, \xi_2) | & 0 \leq x_1 \leq 2, 0 \leq x_2 \leq 1, 0.05 \leq \xi_1 \leq 0.45, \\
    & 0.5 \leq \xi_2 \leq 2.5, (x_1-1.5)^2+(x_2-0.5)^2 \geq \xi_1^2 \}.
    \end{aligned}
\end{equation*} 

Following \cite{yin2024aonn}, we use the necessary conditions for the minimizer of the PDE system described in Section~\ref{section_setup} to find the optimal solution of the parametric OCP. Specifically, the solution of \eqref{eq:test1} is obtained by solving the following KKT systems:
\begin{equation}\label{eq:test1KKT}
    \left\{\begin{aligned}
        &\left\{\begin{aligned}
        - \Delta y(\bx, \bxi) &= u(\bx, \bxi)  &&\text{in} \, \Omega, \\
        y(\bx, \bxi) &= 1   &&\text{on} \, \partial \Omega, 
        \end{aligned}\right.\\
        &\left\{\begin{aligned}
        - \Delta p(\bx, \bxi) &= y(\bx, \bxi) - y_d(\bx, \bxi)   &&\text{in} \, \Omega, \\
        p(\bx, \bxi) &= 0  &&\text{on} \, \partial \Omega, 
        \end{aligned}\right.\\
        &u(\bx, \bxi) = -\frac{1}{\alpha}P_{[u_a,u_b]}\big(p(\bx, \bxi)\big) \quad \text{in} \, \Omega, \\
        \end{aligned}\right. \\
\end{equation}
where $p(\bx, \bxi)$ is the adjoint variable and $P_{[u_a,u_b]}\big(p(\bx, \bxi)\big)$ is defined as:
\begin{equation*}
        P_{[u_a,u_b]}\big(p(\bx, \bxi)\big) = \left\{\begin{aligned}
        &u_b,  &&\text{if} \;\, p(\bx, \bxi)) > u_b, \\
        &p(\bx, \bxi),  &&\text{if} \;\,u_a \leq p(\bx, \bxi)) \leq u_b, \\
        &u_a,  &&\text{if} \;\, p(\bx, \bxi)) < u_a. \\
         \end{aligned}\right. \\
\end{equation*}

In this problem, we employ penalty-free techniques to ensure that the neural networks naturally satisfy the homogeneous Dirichlet boundary conditions in both the state equations and the adjoint equations. Specifically, we first introduce three neural networks $\hat{y}_I(\bx,\bxi; \theta_{y_I}), \hat{u}(\bx,\bxi; \theta_u),$ and $\hat{p}_I(\bx,\bxi; \theta_{p_I})$, each consisting of 6 fully connected layers with 25 neurons in each hidden layer. And following~\cite{yin2024aonn}, we then define a length factor function as:
\begin{equation*}
    l(\bx, \bxi) = x_1(2-x_1)x_2(1-x_2)\big(\xi_1^2-(x_1-1.5)^2-(x_2-0.5)^2\big).
\end{equation*}
Using these neural networks and the length factor function, the state function and the adjoint function can be then approximated as:
\begin{equation*}
    \hat{y}(\bx,\bxi; \theta_y) \approx l(\bx,\bxi)\hat{y}_I(\bx,\bxi; \theta_{y_I})+1, \quad \hat{p}(\bx,\bxi; \theta_p) \approx l(\bx,\bxi)\hat{p}_I(\bx,\bxi; \theta_{p_I}).
\end{equation*}
Note that both $\hat{y}(\bx,\bxi; \theta_y)$ and $\hat{p}(\bx,\bxi; \theta_p)$ satisfy the Dirichlet boundary conditions of \eqref{eq:test1KKT}. As a result, the training difficulties associated with enforcing boundary conditions are eliminated, thereby improving accuracy and robustness for complex geometries.

The relevant parameter settings in Algorithm~\ref{alg:DAS-AONN} for this numerical example are set as follows: the decay factor $\gamma$ is set to $0.985$; the initial step size $c^{(0)}$ is set to $100$; the initial number of epochs $n^{(0)}$ is set to $200$; the maximum number of epochs $N_{\text{ep}}$ is set to $2000$; the batch size $m$ is set to $2000$; the number of adaptive iterations $N_{\text{adaptive}}$ is set to $10$; and $n_{\text{aug}}$ is set to $4$ to ensure the convergence of the algorithm. For KRnet, we set $K = 2$ and take $L = 6$ affine coupling layers. For each affine coupling layer, a two-layer fully connected neural network is used, with 24 neurons in each hidden layer. KRnet is trained by the ADAM optimizer with a learning rate 0.0001, the number of epochs for training $p_{\text{KRnet}}(x,\xi,\theta_f)$ is set to $N_{\text{ep}} = 2000$.

The collocation points in the initial training set $\mathrm{S}_{\Omega, 0}$ are sampled uniformly, with $n_r = \vert \mathrm{S}_{\Omega, 0} \vert$ used during the adaptive sampling procedure. In addition, all three neural networks $\hat{y}(\bx,\bxi; \theta_y), \hat{u}(\bx,\bxi; \theta_u),$ and $\hat{p}(\bx,\bxi; \theta_p)$ are trained using the BFGS method~\cite{mcfall2009artificial}.

To validate the performance of the adaptive AONN method, we conduct a comparative analysis with $\text{DAS}^2$ and AONN. The configurations of the neural networks for the $\text{DAS}^2$ and AONN methods are the same as those of the adaptive AONN. We adopt the dolfin-adjoint framework \cite{mitusch2019dolfin} to solve the optimal control problem with fixed parameters, treating it as the ground truth. The finite element solutions are evaluated on a 256 $\times$ 128 grid for the physical domain and on an $11 \times 11$ grid for the parametric domain corresponding to $\bxi$.

\begin{figure}[htbp]
\centering
\includegraphics[width=1\linewidth,trim=50 0 50 0,clip]{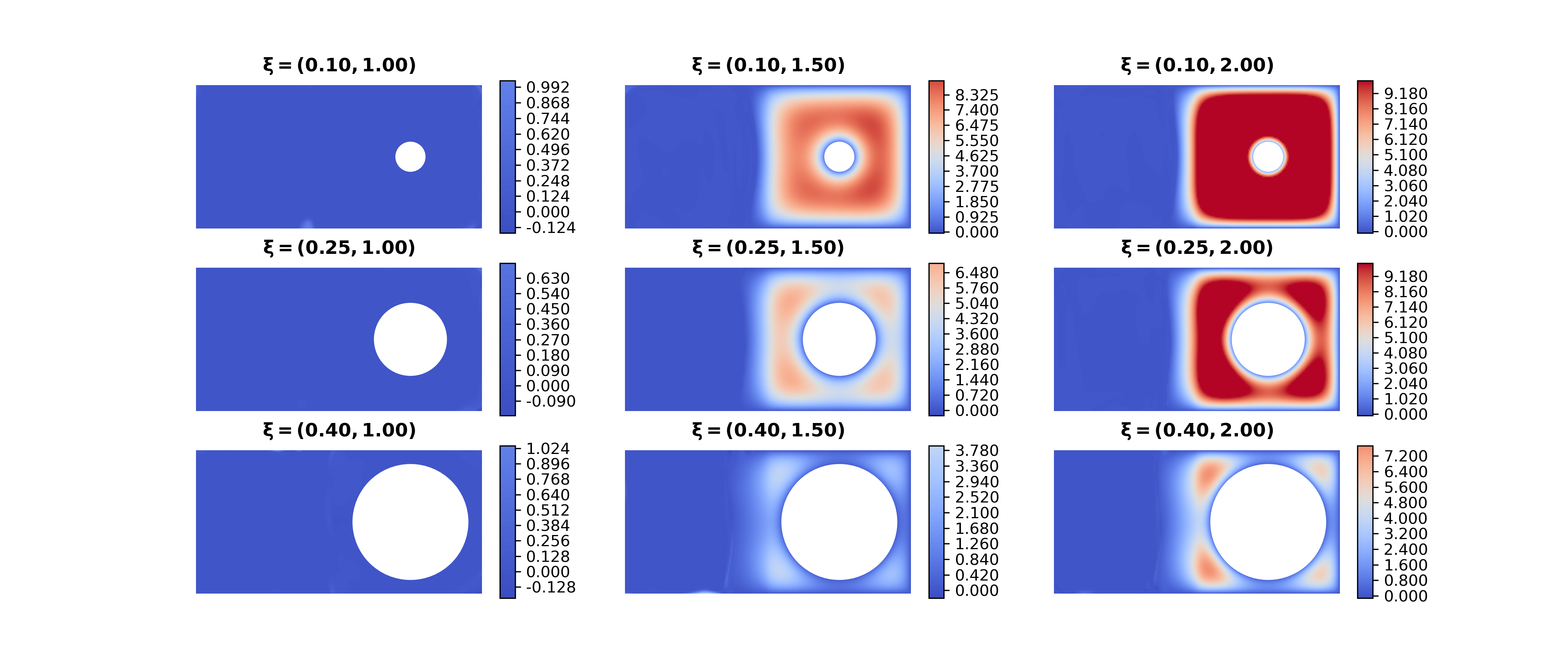}
\caption{adaptive AONN results for the control $u(\bx,\bxi)$ corresponding to different $\bxi = (\xi_1,\xi_2)$ in Test 1.}
\label{DASAONNdiffmu}
\end{figure} 

In Figure~\ref{DASAONNdiffmu}, we show the optimal control solution obtained using adaptive AONN under various parameter settings, respectively. As shown in the figure, when the geometric parameter $\xi_1$ increases, the control generally decreases. In contrast, it increases as the desired state parameter $\xi_2$ increases. This behavior matches practical scenarios: higher temperatures can be applied in smaller tumor regions, while a higher heat source intensity is needed when targeting elevated temperature fields. Especially when $\xi_2 = 1$, indicating that the target temperature is uniform throughout the region, the results from adaptive AONN demonstrate precise agreement with theoretical expectations, as shown in the left columns of the figures.

\begin{figure}[htbp]
        \centering
        \setlength{\tabcolsep}{0pt}
        \renewcommand{\arraystretch}{0}
        
        \begin{tabular}{
                >{\centering\arraybackslash}m{2cm}  
                >{\centering\arraybackslash}m{0.43\textwidth}  
                >{\centering\arraybackslash}m{0.43\textwidth}  
            }
            
            & \text{Control function} & \text{Absolute error} \\[4pt]
            
            \text{$\text{DAS}^2$} &
            \includegraphics[width=0.43\textwidth]{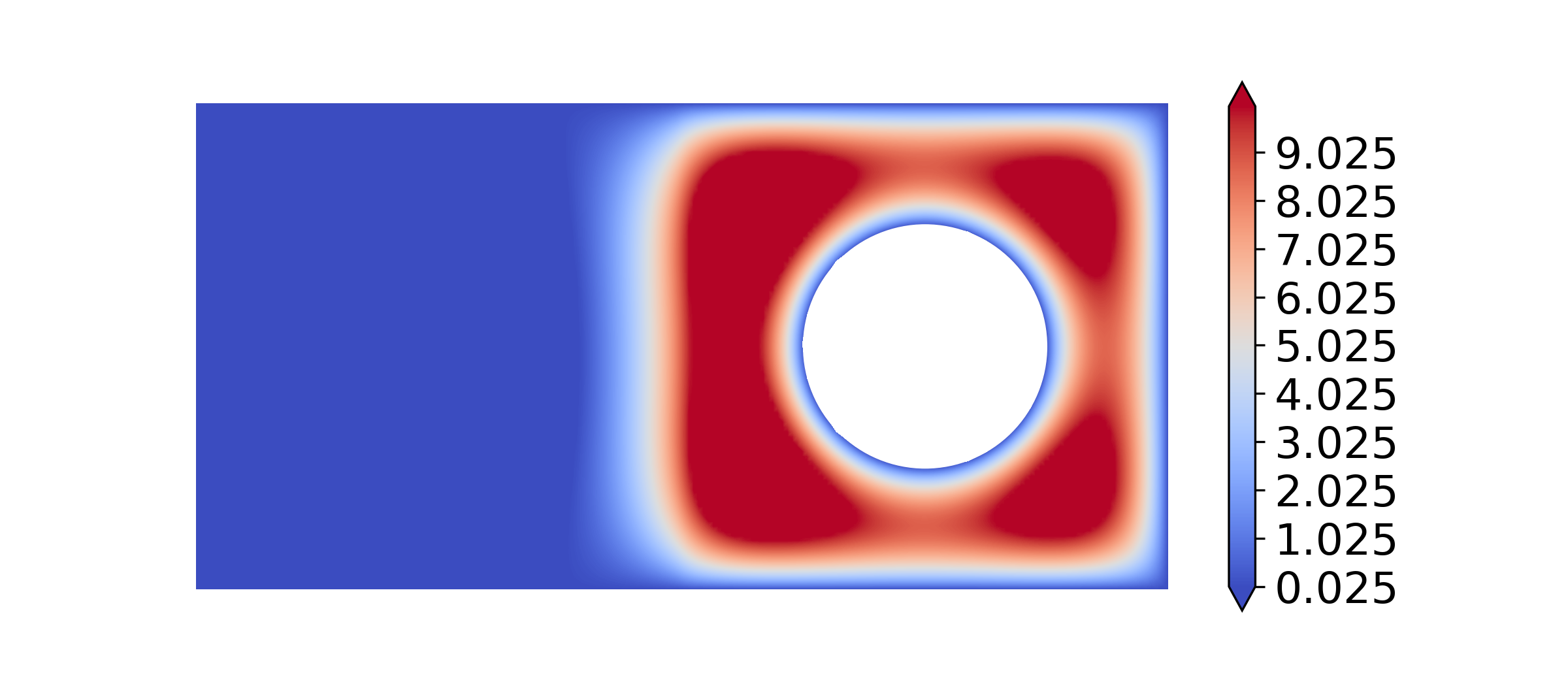} &
            \includegraphics[width=0.43\textwidth]{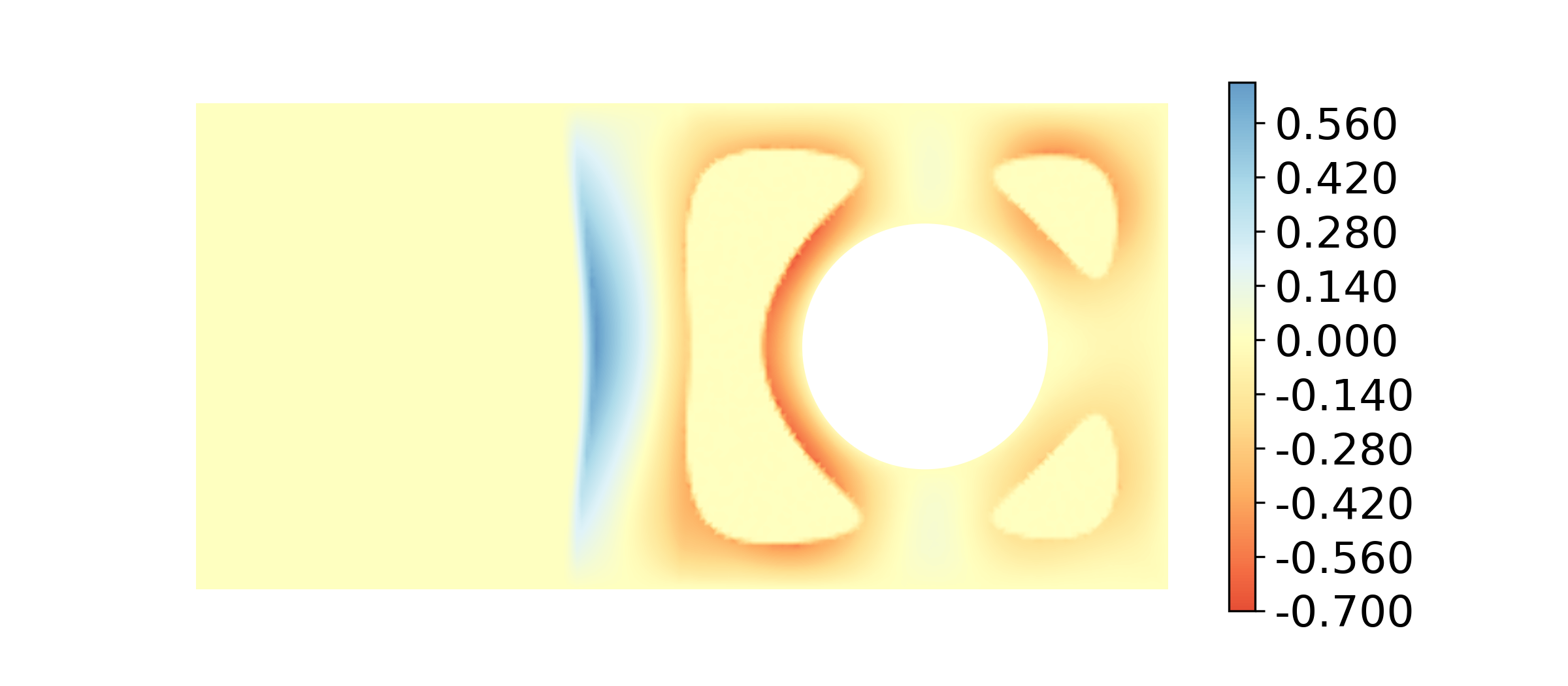} \\[6pt]
            
            \text{AONN} &
            \includegraphics[width=0.43\textwidth]{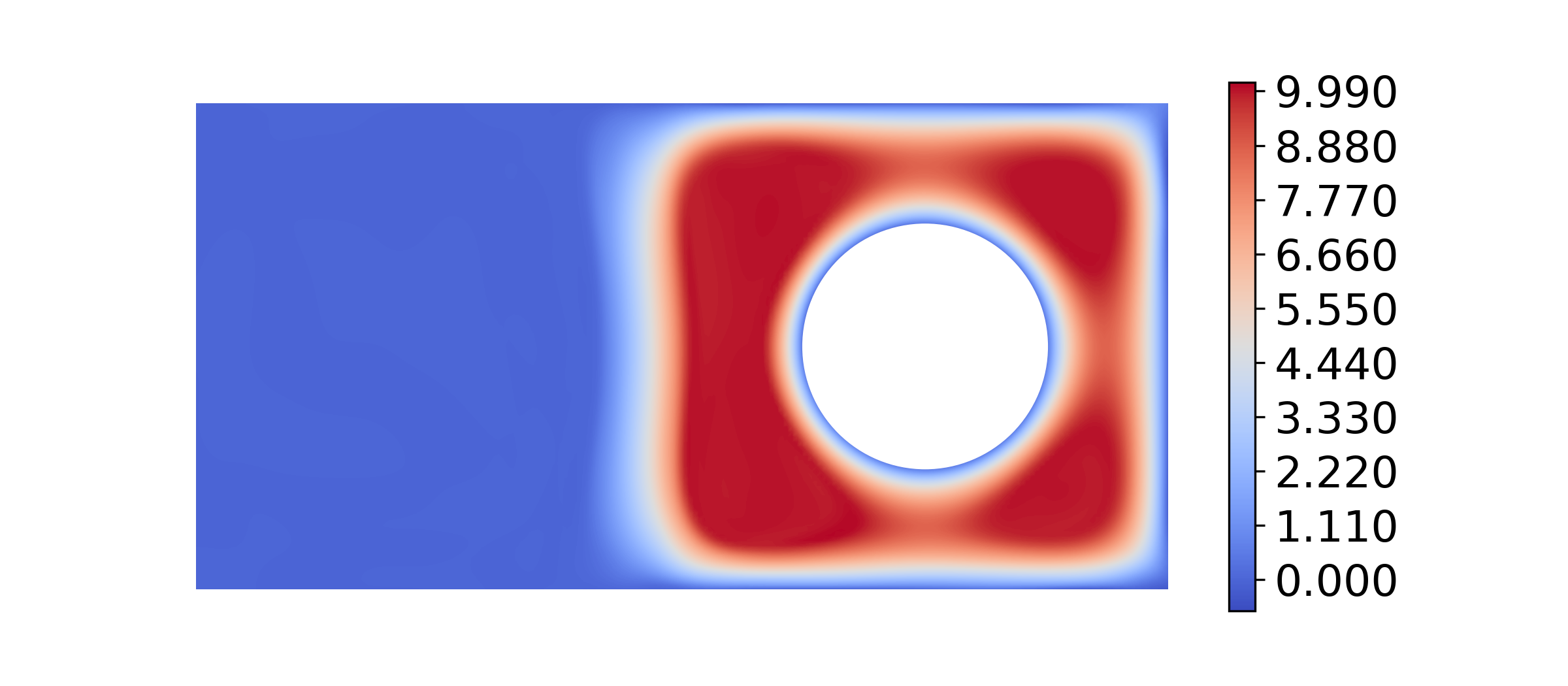} &
            \includegraphics[width=0.43\textwidth]{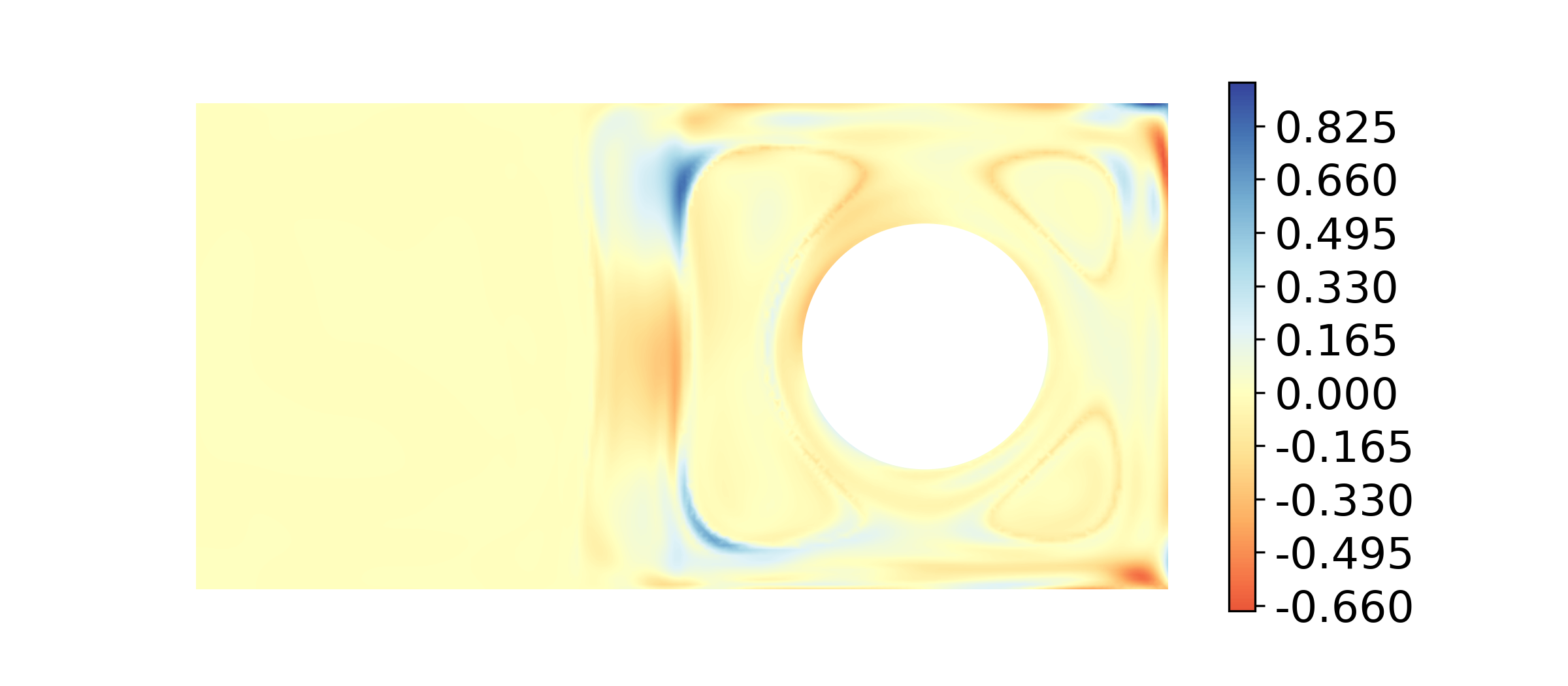} \\[6pt]
            
            \makecell{Adaptive\\AONN} &
            \includegraphics[width=0.43\textwidth]{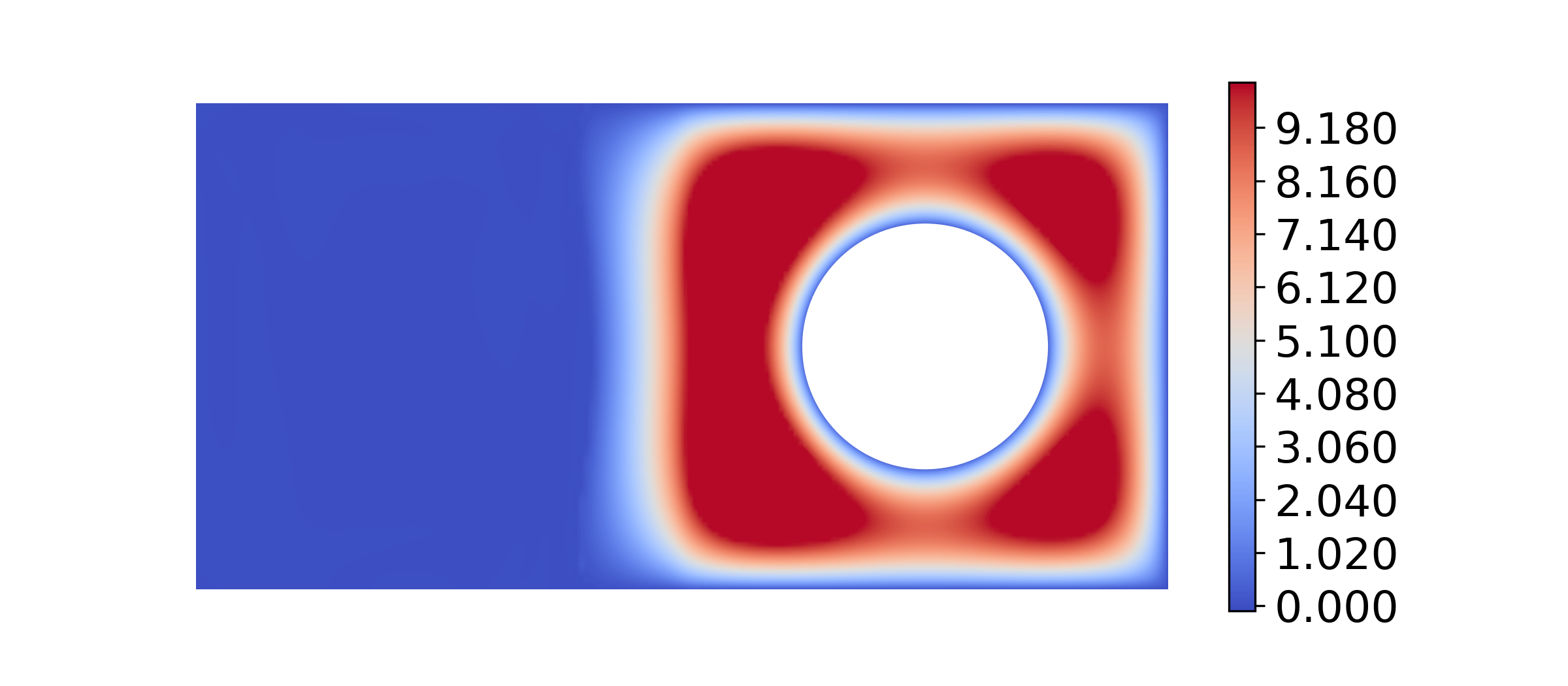} &
            \includegraphics[width=0.43\textwidth]{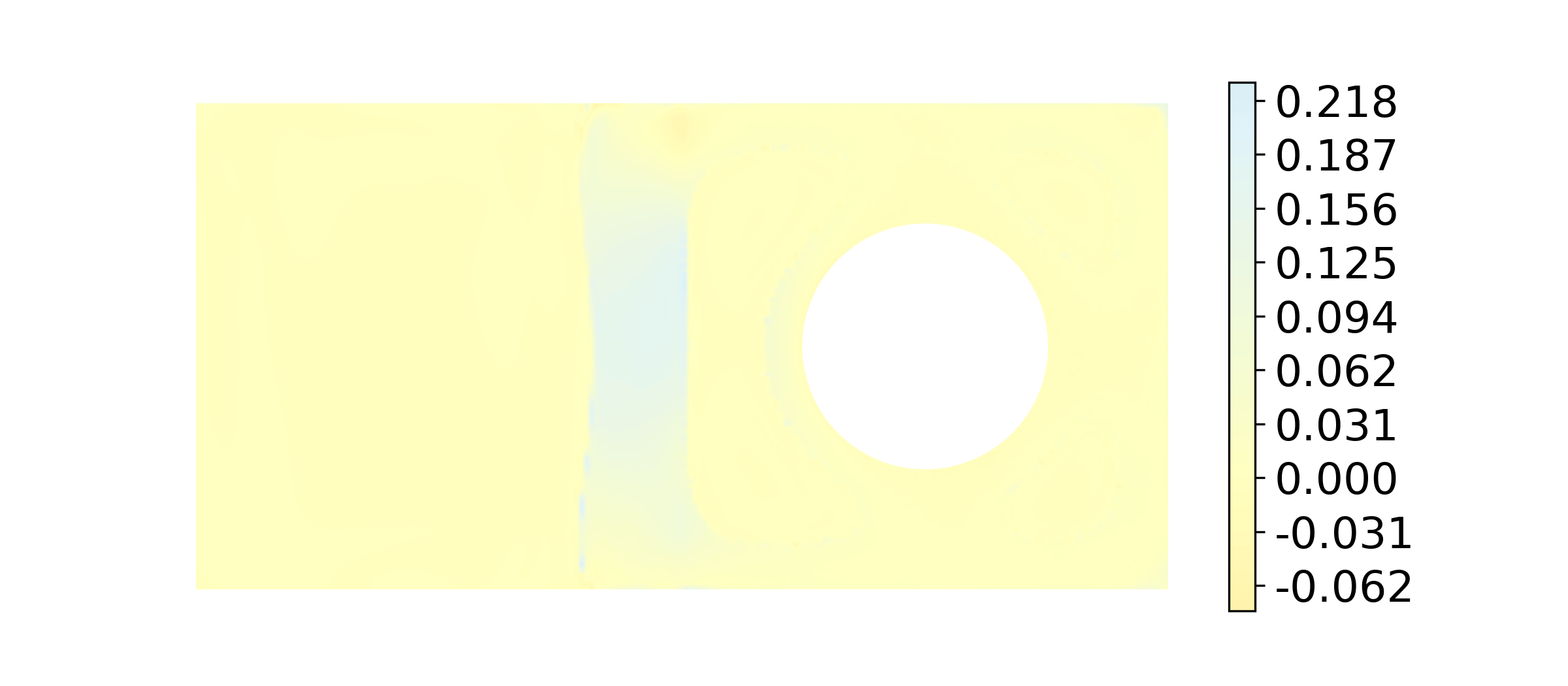}
        \end{tabular}
        
        \caption{Results for the control $u(\bx,\bxi)$ and its absolute error for the fixed parameter $\bxi = (0.25,2.00)$ in Test 1. }
    \label{fixedmu}
    \end{figure}

Figure~\ref{fixedmu} presents a comparison of $\text{DAS}^2$, AONN, and adaptive AONN for the parameter setting $\bxi = (0.25, 2.00)$. The left column displays the results of the control variable $u$ obtained by the three methods. It can be seen from the figure that all three methods have achieved similar distributions, although their value ranges differ. This implies that $\text{DAS}^2$, AONN, and adaptive AONN are all capable of generating approximate solutions, but they differ in precision. The right column presents the pointwise errors between the solutions obtained by the three methods and the dolfin-adjoint reference. It is clear that the adaptive AONN solution achieves the highest accuracy, with a relative $l_2$ error of $0.006$ compared to the reference solution, while the relative $l_2$ error for $\text{DAS}^2$ and AONN are $0.025$ and $0.024$, respectively.

\begin{figure}[htbp]
	\begin{center}
		\subfloat[Relative $l_2$ error w.r.t. sample size.]{
			\includegraphics[width=0.45\linewidth]{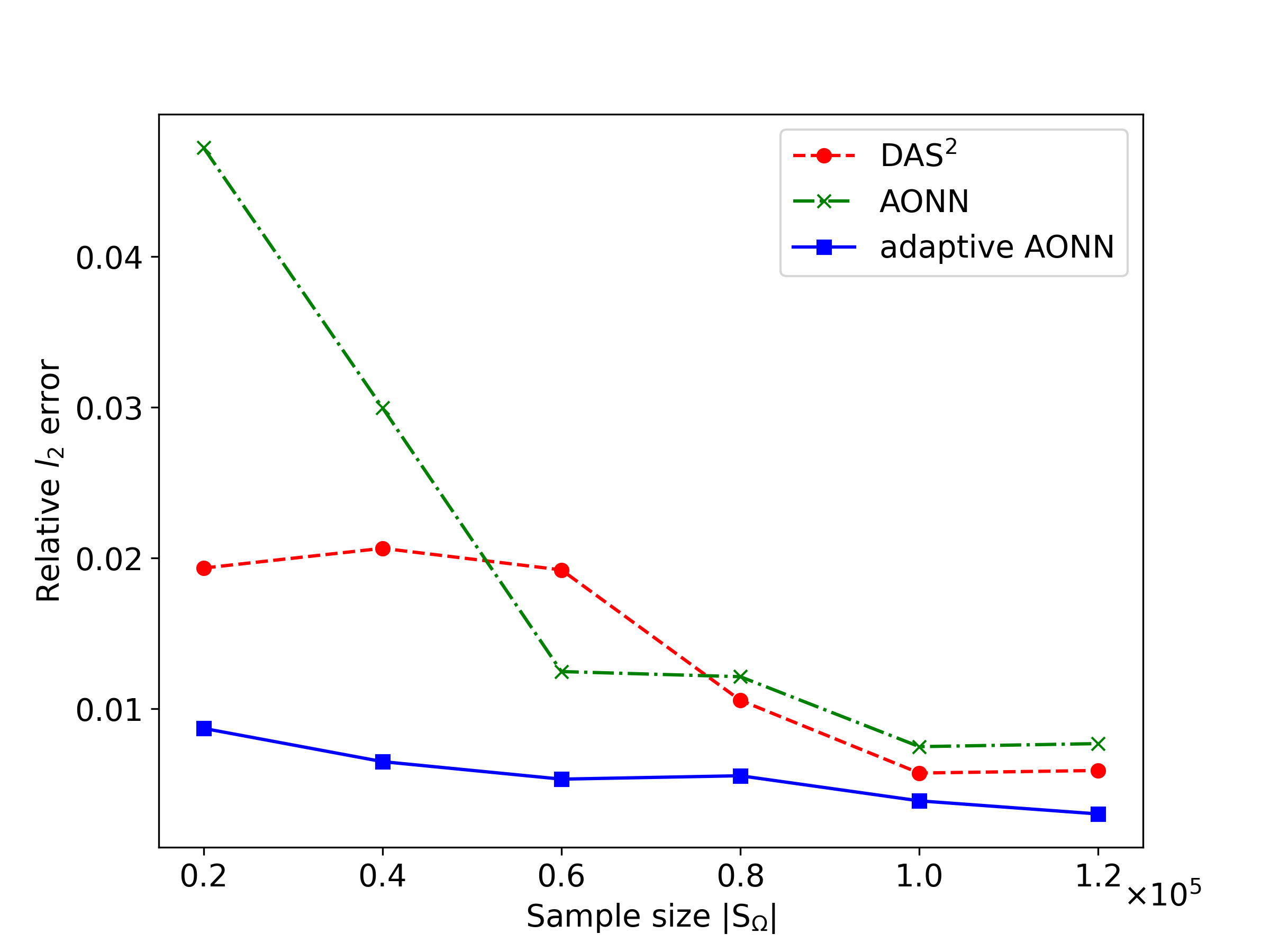}
            \label{5.1}
		}
		\subfloat[Relative $l_2$ error w.r.t. the number of training epochs.]{
			\includegraphics[width=0.45\linewidth]{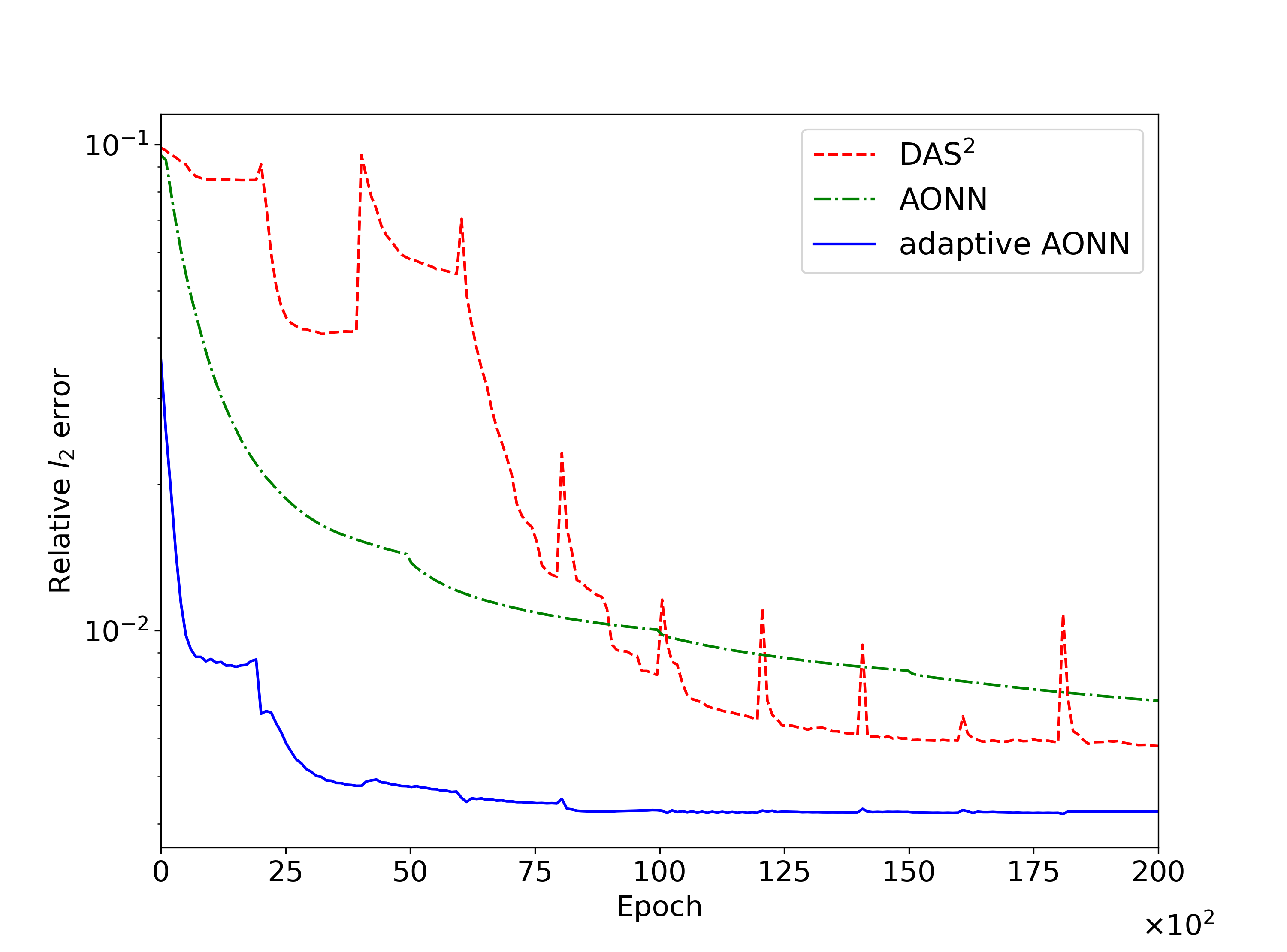}
            \label{5.2}
		}
	\end{center}
    \caption{Relative $l_2$ error w.r.t. sample size and the number of training epochs in Test 1.}\label{reltive_error_sample_fig}
\end{figure}

Figure \ref{5.1} plots the relative $l_2$ errors of the solutions generated by different training methods across varying numbers of collocation points. For adaptive AONN and $\text{DAS}^2$, the initial training set size $\vert \mathrm{S}_{\Omega, 0} \vert = n_r$ is set to  $2\times10^3, 4\times10^3, 6\times10^3, 8\times10^3, 1\times10^4$, and $1.2\times10^4$, corresponding to total training sizes of  $2\times10^4, 4\times10^4, 6\times10^4, 8\times10^4, 1\times10^5$, and $1.2\times10^5$, respectively. For AONN, the number of training points is chosen to match these settings for a fair comparison. The curves show that adaptive AONN outperforms both $\text{DAS}^2$ and AONN, particularly in low-data regimes: it achieves comparable accuracy with sparse samples whereas $\text{DAS}^2$ and AONN require substantial more data ($\vert \mathrm{S}_{\Omega} \vert = 1.0\times10^5$) to reach a similar precision. This highlights the data efficiency and robustness of adaptive AONN in sample-limited scenarios.

Figure~\ref{5.2} shows the error evolution of the three methods during their adaptive iteration processes, with the initial training set size set to $n_r = 1 \times 10^4$. From the figure, we can see that AONN exhibits a steady decrease in error as the number of epochs increases, whereas $\text{DAS}^2$ and adaptive AONN display noticeable spikes at the beginning of each adaptive stage. This behavior arises from modifications to the sampling set in deep adaptive sampling, since newly added points tend to introduce larger errors. A closer comparison of the two adaptive strategies reveals that $\text{DAS}^2$ reduces the error gradually throughout the entire training process, while adaptive AONN achieves rapid error reduction within the first two adaptive stages and then remains essentially stable. These results indicate that adaptive AONN reaches convergence with fewer adaptive iterations.

\begin{figure}[h]
    \centering
    \includegraphics[width=0.8\textwidth]{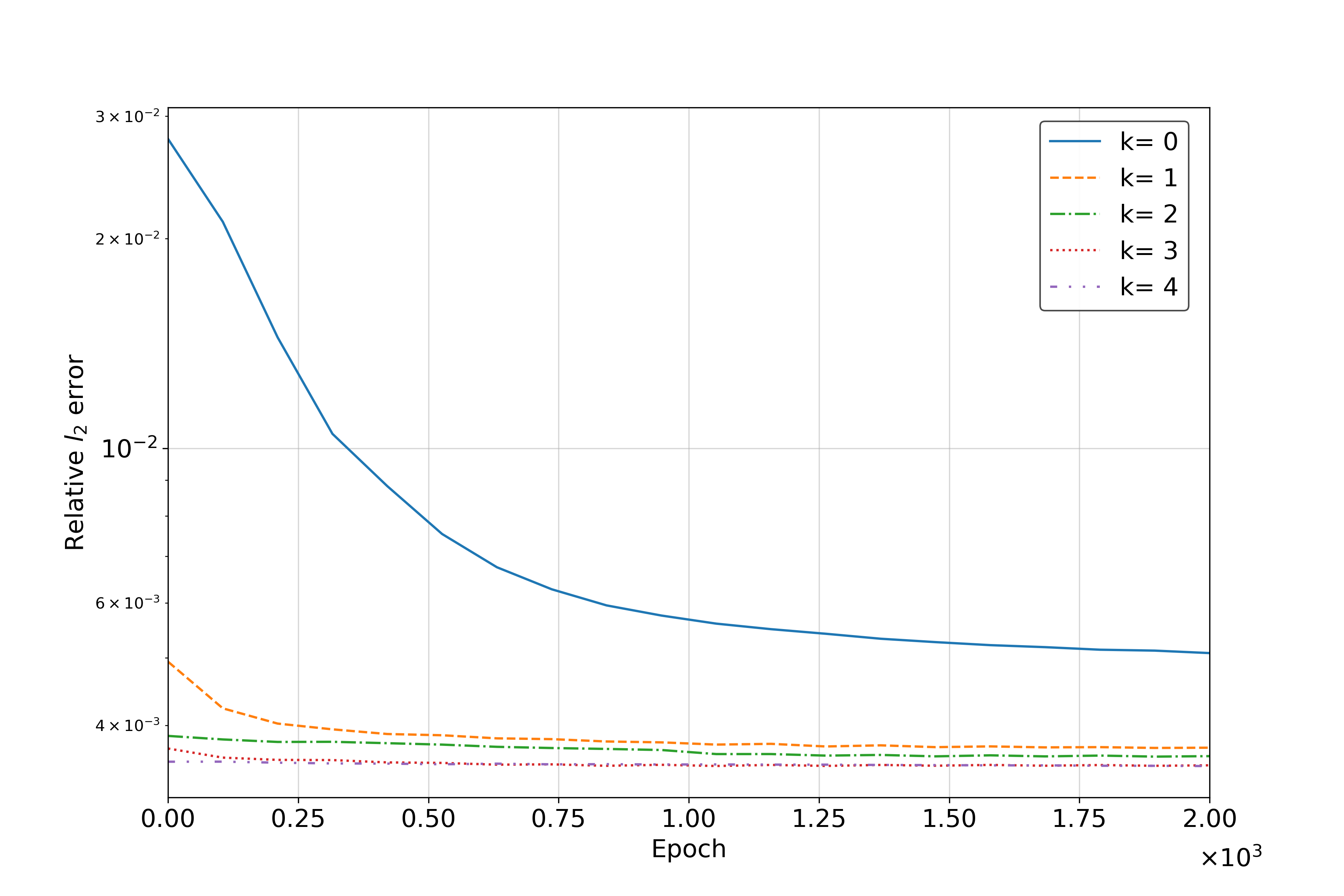}
    \caption{adaptive AONN error curves during different adaptive iterations procedure in Test 1.}
    \label{dasaonn_epoch_error}
\end{figure}

Figure~\ref{dasaonn_epoch_error} further displays the error trajectories of adaptive AONN across different adaptive stages, where substantial error reduction occurs only in the first two iterations. A detailed analysis of this behavior reveals the following mechanism: during the initial iteration, both $\text{DAS}^2$ and adaptive AONN employ randomly generated collocation points for training. At this stage, the difference in surrogate model structure determines performance. As shown in Figure~\ref{reltive_error_sample_fig} (b), adaptive AONN significantly outperforms $\text{DAS}^2$ due to its enhanced surrogate modeling through the AONN architecture, compared with the standard fully-connected neural network used in $\text{DAS}^2$. In the second iteration, adaptive AONN exhibits another pronounced error reduction because the deep generative model effectively samples new collocation points from high-residual regions identified in the first iteration, which makes a significant contribution to surrogate model training. After these two iterations, adaptive AONN nearly reaches numerical convergence, whereas $\text{DAS}^2$ requires more adaptive iterations to reach comparable accuracy.

\subsection{Test 2: Optimal control for the Stokes equation with physical parametrization}
In this problem, we consider the following parametric OCP:
\begin{equation*}
    \min_{g(\bx,\xi),u(\bx,\xi),p(\bx,\xi)} J\big(g(\bx,\xi),u(\bx,\xi), p(\bx,\xi)\big) = \frac{1}{2}\int_{\Omega}\nabla u(\bx,\xi) \cdot \nabla u(\bx,\xi) \,\mathrm{d}x + \frac{\alpha}{2}\int_{\partial\Omega_{\text{circle}}} (g(\bx,\xi))^2 \, \mathrm{d}s,
\end{equation*}
subject to the following Stokes equations:
\begin{equation}\label{test:stokes_state}
\left\{
    \begin{aligned}
        -\xi \Delta u(\bx,\xi) + {\nabla} p(\bx,\xi) &= 0 \qquad &\text{in}& \; \Omega, \\
        \text{div} \,u(\bx,\xi) &= 0 \qquad &\text{in}& \; \Omega,\\
        u(\bx,\xi) &= g(\bx,\xi)  &\text{on}& \; \partial \Omega_{\text{circle}}, \\
        u(\bx,\xi) &= f(\bx) \ & \text{on}& \; \partial \Omega_{\text{in}}, \\
        u(\bx,\xi) &= 0 \ & \text{on}& \; \partial \Omega_{\text{walls}}, \\
        p(\bx,\xi) &= 0 \ & \text{on}& \; \partial \Omega_{\text{out}}, \\
    \end{aligned}
\right.
\end{equation}

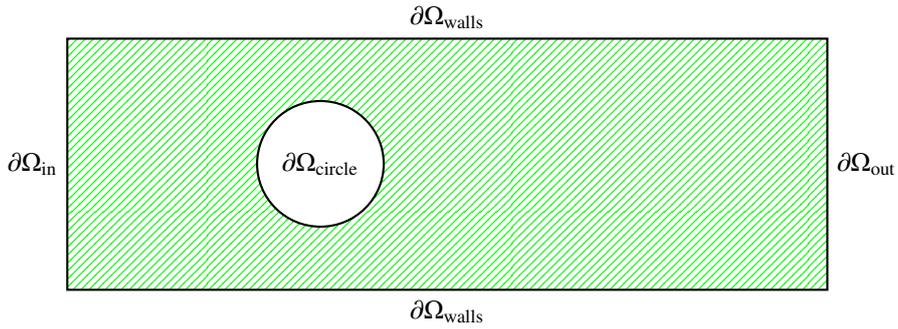
\begin{figure}[htbp]
\centering
\begin{tikzpicture}[thick]
  \draw[pattern=north east lines, pattern color=green, even odd rule] (0, 0) rectangle (10, 3.33);
  \fill[white!20] (3.33,1.67) circle (0.833);
  \draw (3.33,1.67) circle (0.833);
  \node at (3.33,1.67) {$\partial \Omega_{\text{circle}}$}; 
  \node[left, align=center] at (0,1.67) {$\partial \Omega_{\text{in}}$};
  \node[right, align=center] at (10,1.67) {$\partial \Omega_{\text{out}}$};
  \node[above] at (5,3.33) {$\partial \Omega_{\text{walls}}$};
  \node[below] at (5,0) {$\partial \Omega_{\text{walls}}$};
\end{tikzpicture}
\caption{The computational domain $\Omega$ of Test 2.}\label{test:stokes_fig1}
\label{case2domain}
\end{figure}

\noindent where the domain of interest is defined as
\begin{equation*}
    \begin{aligned}
    \Omega \coloneqq \{ (x_1, x_2)\ |\ 0 \leq x_1 \leq 30, 0 \leq x_2 \leq 10, (x_1-10)^2+(x_2-5)^2 \geq 2.5^2 \},
    \end{aligned}
\end{equation*}
and it is visualized in Figure~\ref{test:stokes_fig1}. The parameter $\xi$ represents the reciprocal of the Reynolds number and ranges from $10$ to $1000$. In this problem, the regularization parameter $\alpha$ is set to 10.
The velocity field $u(\bx, \cdot): \Omega \rightarrow \mathbb{R}^2$ and the pressure field $p(\bx, \cdot): \Omega \rightarrow \mathbb{R}$ are both unknowns, while $g(\bx, \cdot)$ denotes the unknown control variable that specifies the Dirichlet boundary condition on the circular boundary. In addition, we impose the following constraint on $g(\bx, \xi)$:
\begin{equation}
    0 \leq g(\bx, \xi) \leq 1 \quad \text { a.e. on } \partial \Omega_{\text{circle}}. \label{G_ad xi}
\end{equation}
The function $f(\bx)$ in the Dirichlet inflow boundary condition is chosen as
\begin{equation*}
    f(\bx) = \frac{1}{25}\,x_2\times(10-x_2) \quad \text{on} \;\partial\Omega_{\text{in}}.
\end{equation*}

This problem models a flow control optimization problem focused on minimizing the loss of flow energy by regulating the velocity field at the circle boundary. To avoid excessive control solutions, a regularizing penalty term is introduced to the objective functional. In this context, our goal is to construct an efficient surrogate model for solving the parametric optimal control problem. 

The adjoint equation of this problem is given by
\begin{equation}\label{test:stokes_adjoint}
    \left\{\begin{aligned}
                -\xi \Delta \lambda(\bx,\xi) + {\nabla} \mu(\bx,\xi) &= - \Delta u(\bx,\xi) &\text{in} \; &\Omega, \\
                \text{div} \, \lambda(\bx,\xi) &= 0  &\text{in} \; &\Omega, \\
                \lambda(\bx,\xi) &= 0 &\text{on} \; &\Omega_{\text{in}}\cup\Omega_{\text{walls}}, \\
                \alpha g(\bx,\xi) + \lambda(\bx,\xi) \cdot \boldsymbol{n} &= 0 &\text{on} \; &\Omega_{\text{circle}}, \\
                \mu(\bx,\xi) &= 0 &\text{on} \; &\Omega_{\text{out}}.
        \end{aligned}\right. 
\end{equation}
where $\boldsymbol{n}$ denotes the unit outward normal vector along the circular boundary, $\lambda(\bx,\xi)$ denotes the adjoint velocity, and $\mu(\bx,\xi)$ denotes the adjoint pressure.

The state equation \eqref{test:stokes_state} and the adjoint equation \eqref{test:stokes_adjoint}, together with the variational inequality
\begin{equation}\label{eq:test2_var}
(\mathrm{d}_gJ\big(u^*(\bx, \xi), g^*(\bx, \xi); \xi \big), h(\bx, \xi) - g^*(\bx, \xi)) \geq 0, \quad \forall h(\bx, \xi) \in G_{ad}(\xi),
\end{equation}
constitute the KKT system~\eqref{KKT}, where the admissible set $G_{ad}(\xi)$ is given by \warn{~\eqref{G_ad xi}}, and the derivative of the objective functional $J$ with respect to the control variable $g(\bx,\xi)$ is expressed as
$$
\mathrm{d}_gJ( \bxi )=\alpha g(\bx,\xi) + \lambda(\bx,\xi) \cdot \boldsymbol{n}.
$$
It should be noted that the variational inequality~\eqref{eq:test2_var} in the KKT system is already implicitly embedded in the boundary condition of the adjoint equation~\eqref{test:stokes_adjoint}. However, since the derivative of the objective functional $\mathrm{d}_g J(\bxi)$ is required to construct the loss function $J_g$, we explicitly state it here.

To solve this problem, we construct five neural networks: $\hat{u}(\bx,\xi;\theta_u)$, $\hat{p}(\bx,\xi;\theta_p)$, $\hat{\lambda}(\bx,\xi;\theta_{\lambda})$, $\hat{\mu}(\bx,\xi;\theta_{\mu})$, and $\hat{g}(\bx,\xi;\theta_g)$. Specifically, $\hat{u}(\bx,\xi;\theta_u)$ and $\hat{g}(\bx,\xi;\theta_g)$ consist of $5$ fully connected layers with $32$ neurons in each hidden layer, and $\hat{p}(\bx,\xi;\theta_p)$ has $5$ fully connected layers with $16$ neurons in each hidden layer. The network architectures for $\hat{\lambda}(\bx,\xi;\theta_{\lambda})$ and $\hat{\mu}(\bx,\xi;\theta_{\mu})$ are identical to those of $\hat{u}(\bx,\xi;\theta_u)$ and $\hat{p}(\bx,\xi;\theta_p)$, respectively. For KRnet, we set $K = 3$ and take $L = 6$ affine coupling layers. Each affine coupling layer has $2$ fully connected layers with $24$ neurons in each hidden layer. 

The inputs in Algorithm~\ref{alg:DAS-AONN} are set as follows: the decay factor $\gamma$ is set to $1$; the initial step size $c^{(0)}$ is set to $1$; the initial number of epochs $n^{(0)}$ is set to 400; the maximum number of epochs $N_{\text{ep}}$ is set to $4000$; the batch size $m$ is set to $4000$; the number of adaptive iterations $N_{\text{adaptive}}$ is set to $10$; and $n_{\text{aug}}$ is set to $2$ to ensure the convergence of the algorithm.

The collocation points in the initial training set $\mathrm{S}_{\Omega, 0}$ are sampled uniformly, with $n_r = \vert \mathrm{S}_{\Omega, 0} \vert=8\times 10^3$ used during the adaptive sampling procedure. In addition, the five neural networks $\hat{u}(\bx,\xi;\theta_u)$, $\hat{p}(\bx,\xi;\theta_p)$, $\hat{\lambda}(\bx,\xi;\theta_{\lambda})$, $\hat{\mu}(\bx,\xi;\theta_{\mu})$ and $\hat{g}(\bx,\xi;\theta_g)$ are trained using the BFGS method~\cite{mcfall2009artificial}, while KRnet is trained using the ADAM optimizer with a learning rate of $0.0001$. For the boundary term, $2500$ boundary points are uniformly sampled on each edge of the rectangle, and $5000$ boundary points are sampled on the circumference of the circle.

To evaluate the performance of the adaptive AONN method, we conduct a comparative analysis with $\text{DAS}^2$ and AONN. The neural network configurations for $\text{DAS}^2$ and AONN are identical to those used in adaptive AONN. We employ the dolfin-adjoint framework \cite{mitusch2019dolfin} to solve the optimal control problem with fixed parameters, treating the resulting solution as the ground truth. The finite element solutions are evaluated on a $1500\times 500$ uniform grid for the physical domain. In this problem, the control variable $g(\bx, \xi)$ is applied only on the boundary of the circle, making it difficult to visualize. Thus, we present the optimal velocity field $u(\bx,\xi)$ and pressure distribution $p(\bx,\xi)$ instead for comparison.

\begin{figure}[htbp]
        \centering
        \setlength{\tabcolsep}{0pt}
        \renewcommand{\arraystretch}{0}
        
        \begin{tabular}{
                >{\centering\arraybackslash}m{2cm}  
                >{\centering\arraybackslash}m{0.43\textwidth}  
                >{\centering\arraybackslash}m{0.43\textwidth}  
            }
            
            & \text{Velocity field} & \text{Pressure distribution} \\[4pt]

            \text{dolfin-adjoint} &
            \includegraphics[width=0.43\textwidth]{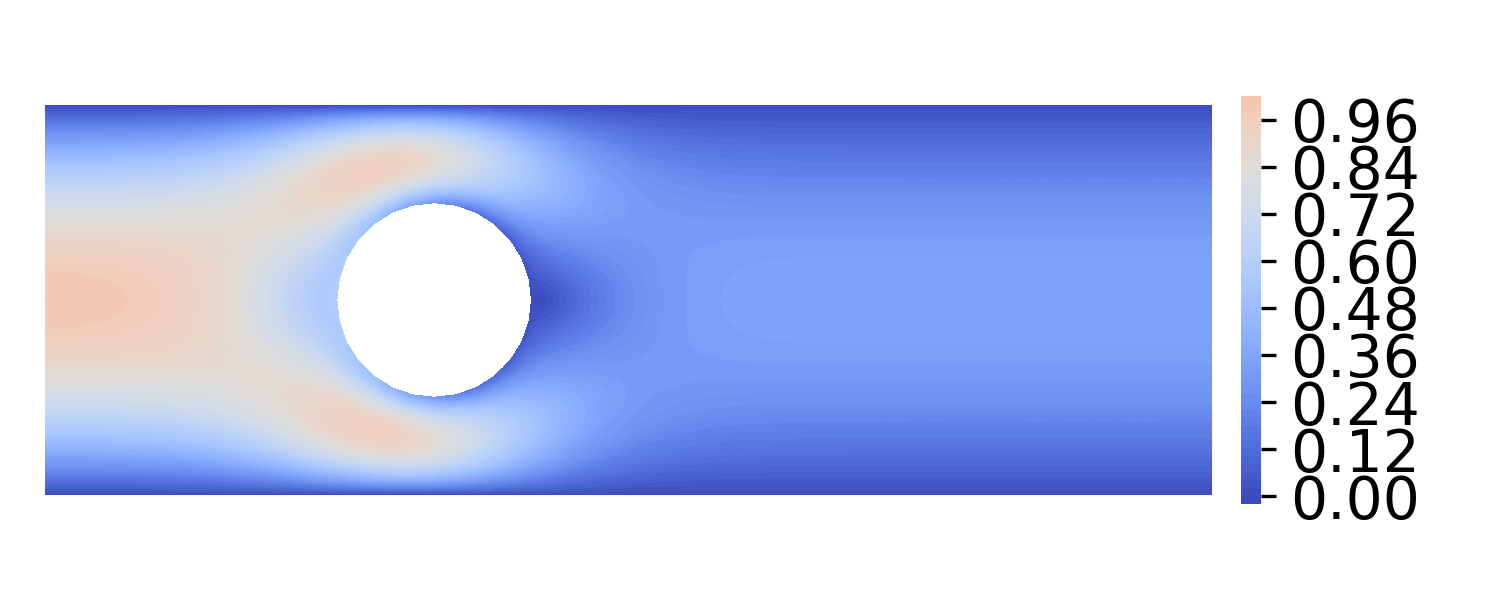} &
            \includegraphics[width=0.43\textwidth]{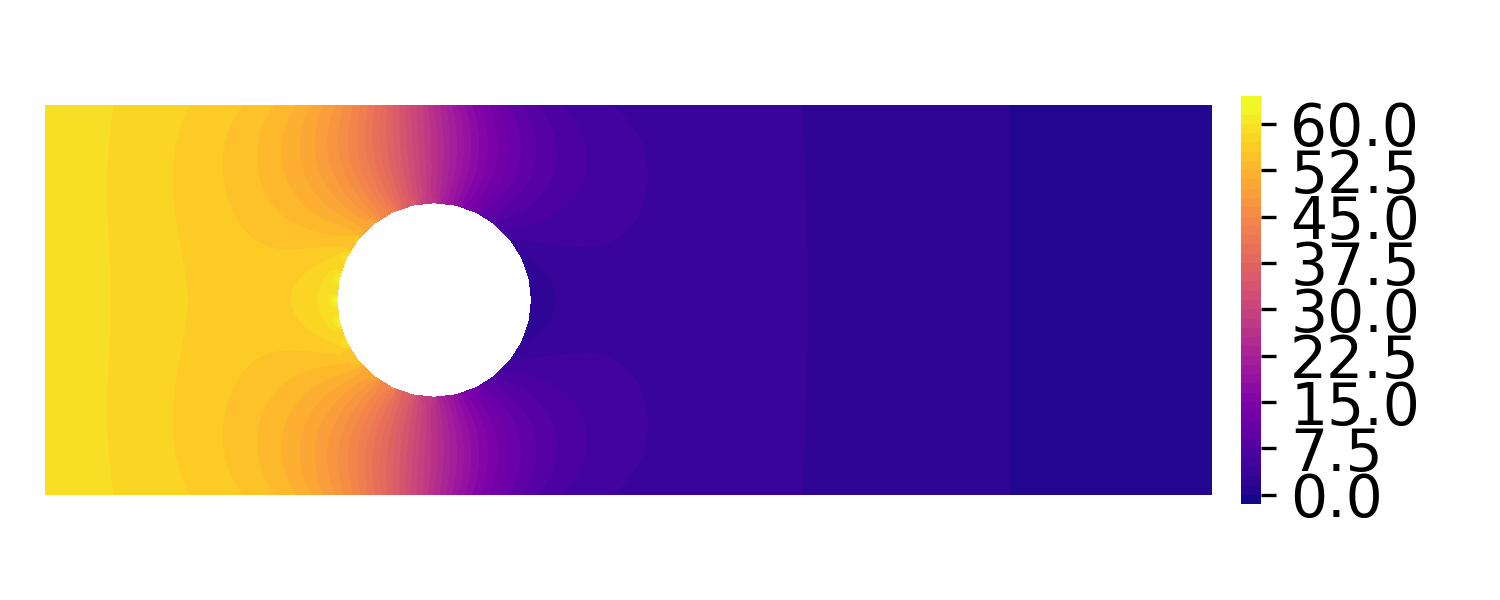} \\[6pt]
            
            \text{$\text{DAS}^2$} &
            \includegraphics[width=0.43\textwidth]{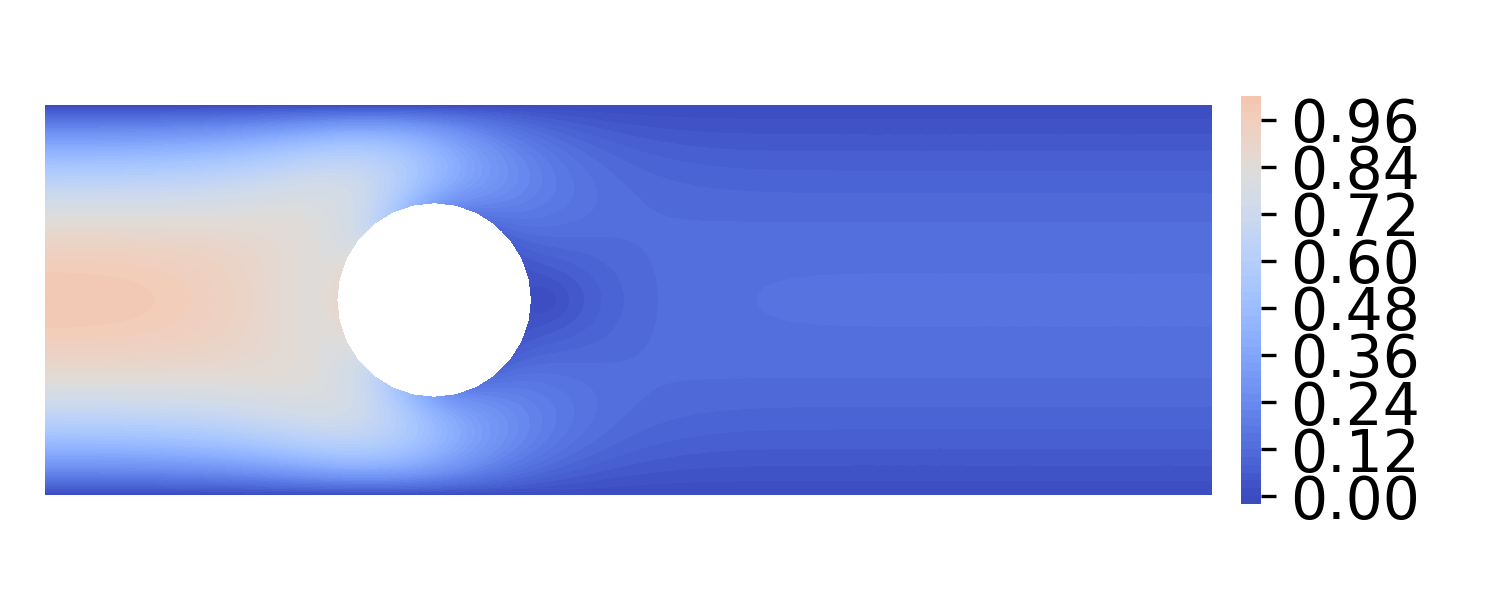} &
            \includegraphics[width=0.43\textwidth]{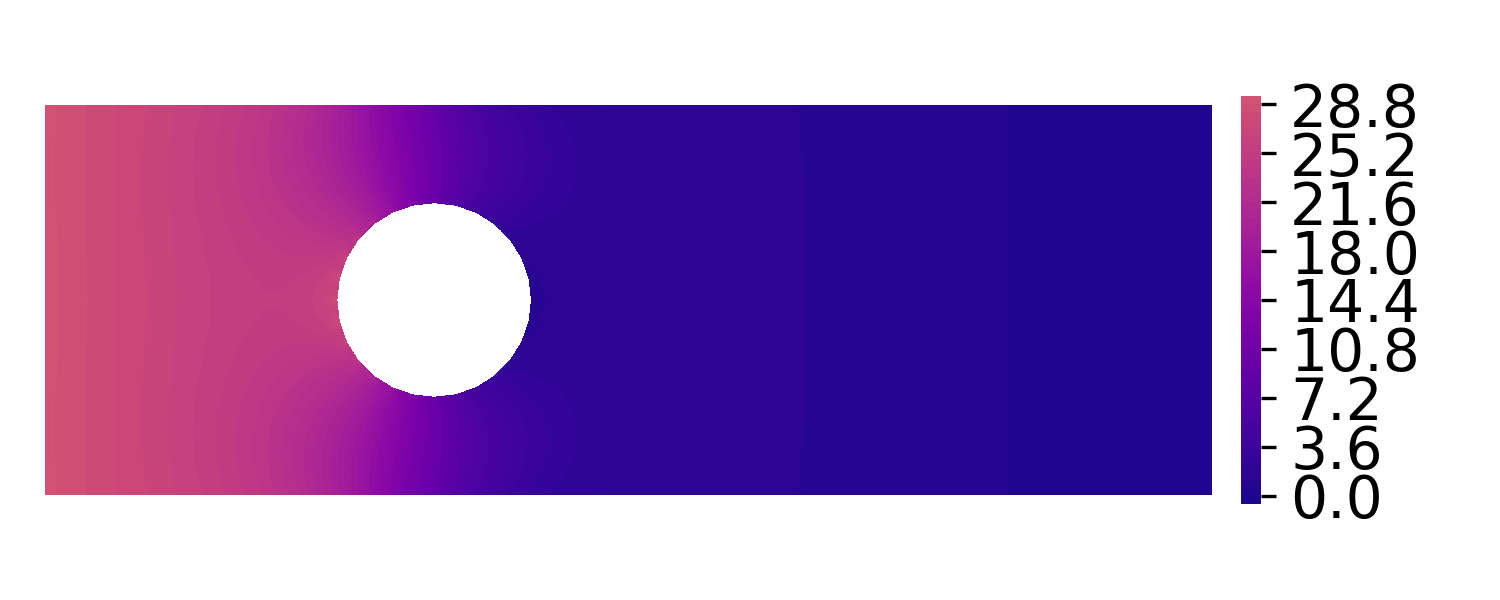} \\[6pt]
            
            \text{AONN} &
            \includegraphics[width=0.43\textwidth]{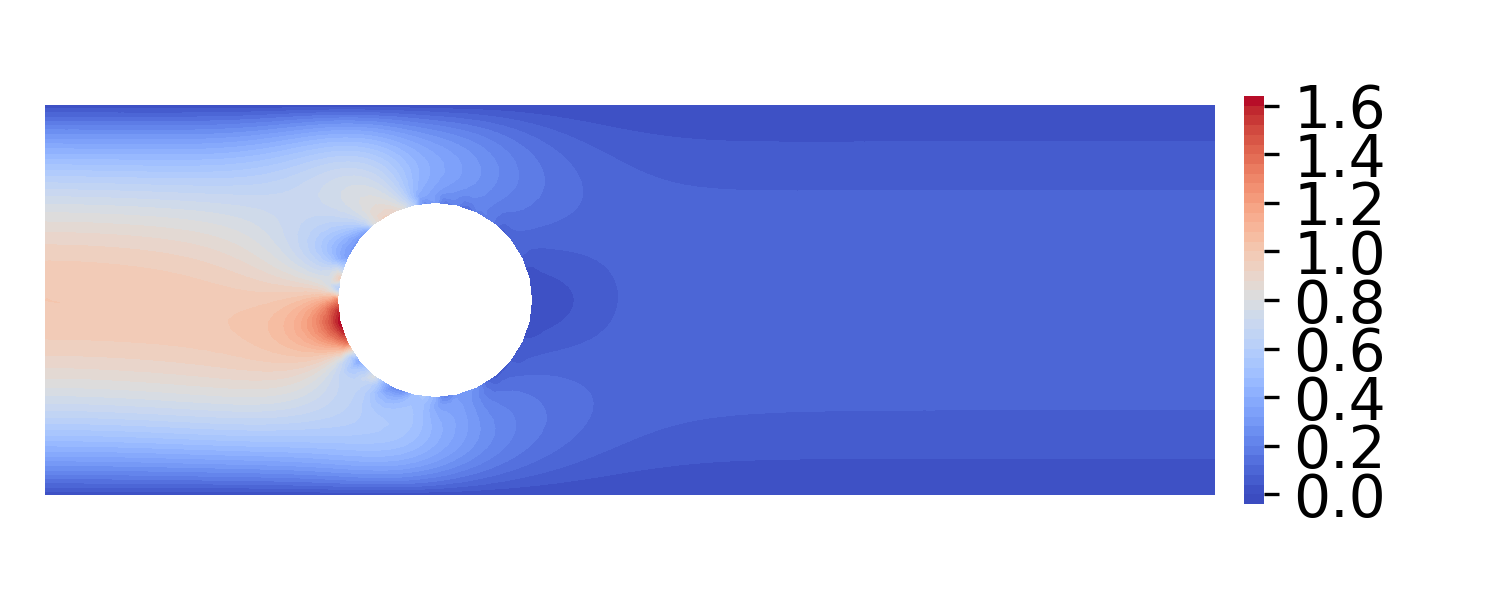} &
            \includegraphics[width=0.43\textwidth]{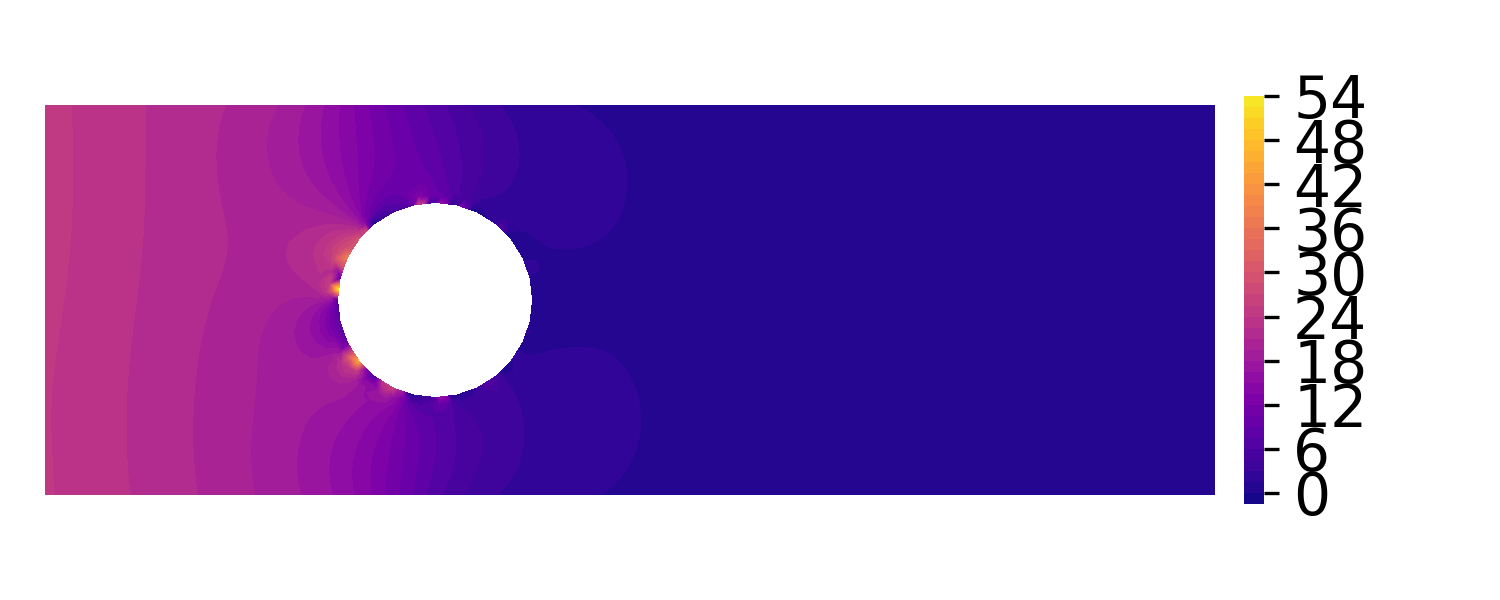} \\[6pt]
            
            \makecell{Adaptive\\ AONN} &
            \includegraphics[width=0.43\textwidth]{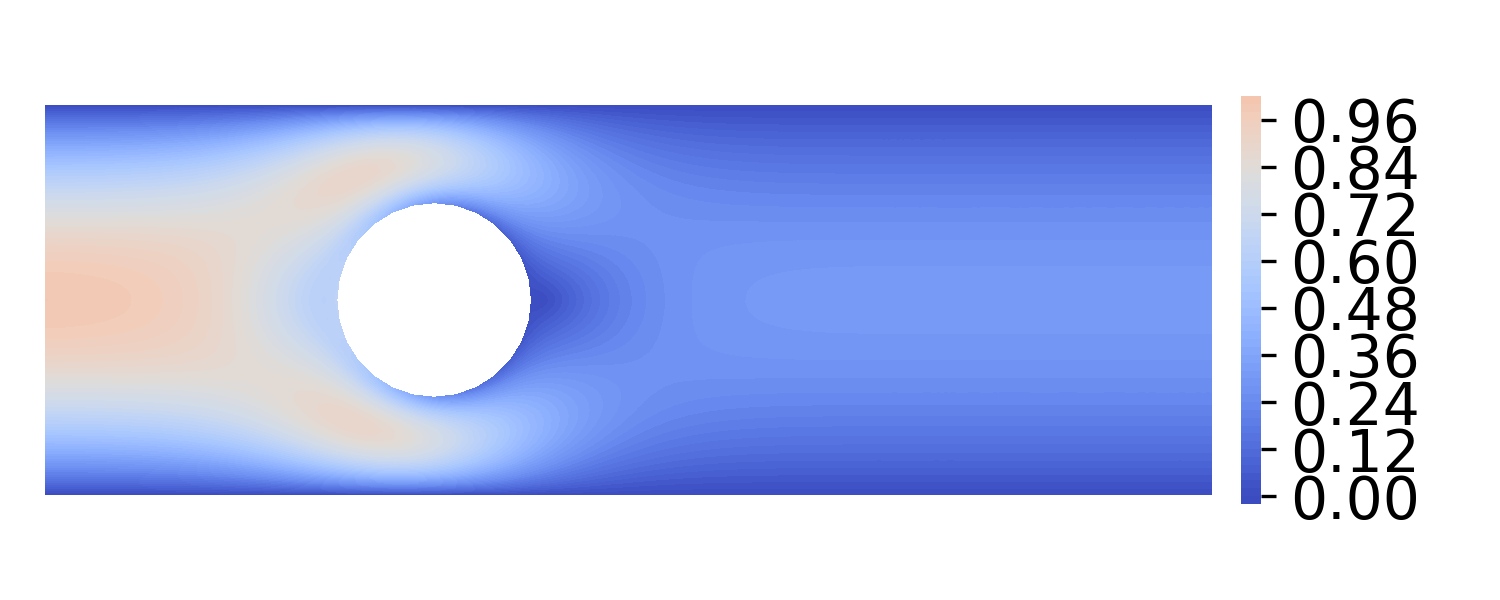} &
            \includegraphics[width=0.43\textwidth]{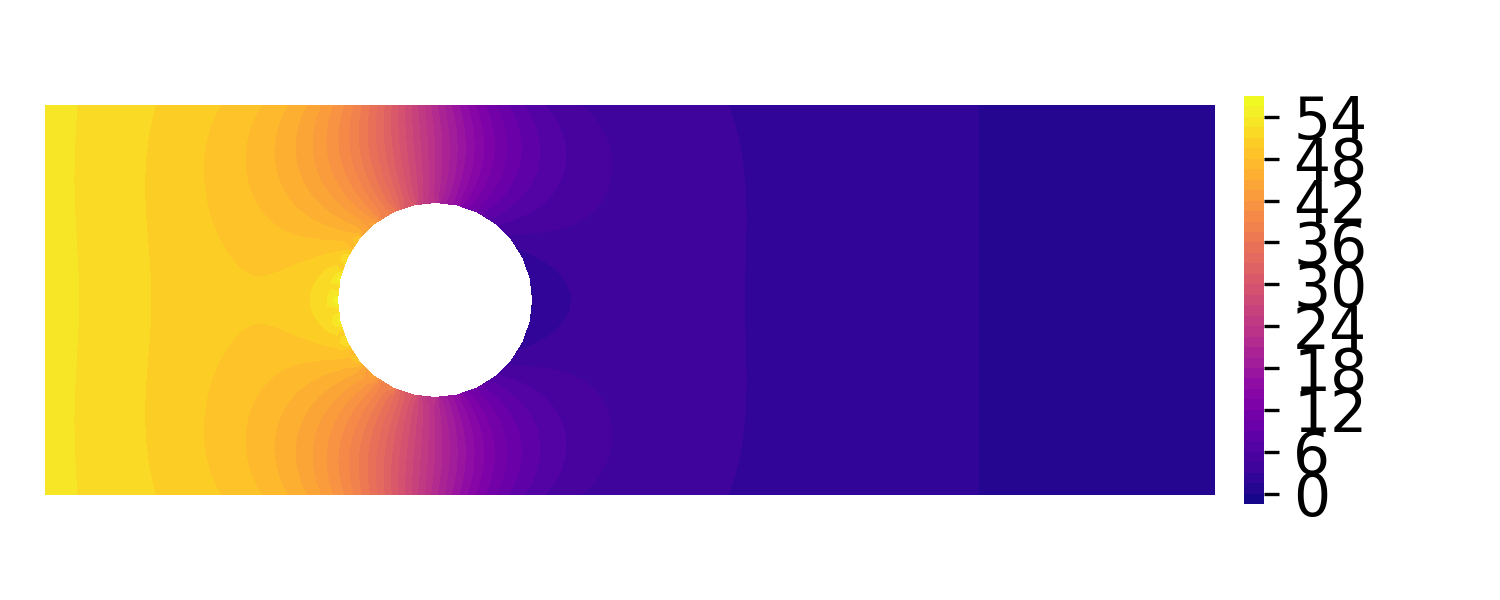}
        \end{tabular}

    \caption{The optimal solutions of the velocity field $u(\bx,\xi)$ and  pressure distribution $p(\bx,\xi)$ for $\xi=10$.} 
    \label{stokes_sol}
\end{figure}

Figure~\ref{stokes_sol} shows the optimal velocity fields $u(\bx,\xi)$ and the optimal pressure distributions $p(\bx,\xi)$ for the parameter value $\xi = 10$. It is evident that both $\text{DAS}^2$ and AONN exhibit pronounced errors in the computation of the optimal velocity and pressure fields compared to the reference solution, particularly as the flow interacts with the left boundary of the circle and evolves downstream. In contrast, adaptive AONN closely matches the finite element solution in both spatial patterns and magnitudes.

\begin{figure}[htbp]
        \centering
        \setlength{\tabcolsep}{0pt}
        \renewcommand{\arraystretch}{0}
        
        \begin{tabular}{
                >{\centering\arraybackslash}m{2cm}  
                >{\centering\arraybackslash}m{0.43\textwidth}  
                >{\centering\arraybackslash}m{0.43\textwidth}  
            }
            
            & \text{Absolute errors of velocity field} & \text{Absolute errors of  pressure distribution} \\[4pt]
            
            \text{$\text{DAS}^2$} &
            \includegraphics[width=0.43\textwidth]{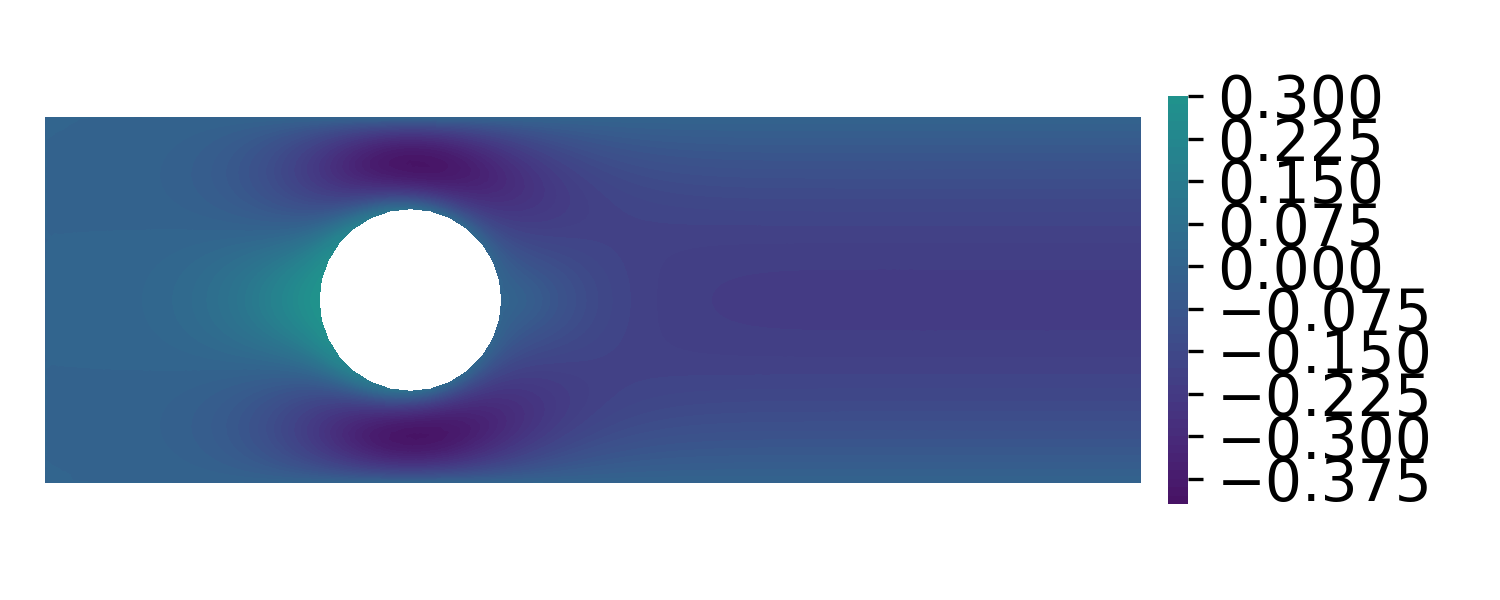} &
            \includegraphics[width=0.43\textwidth]{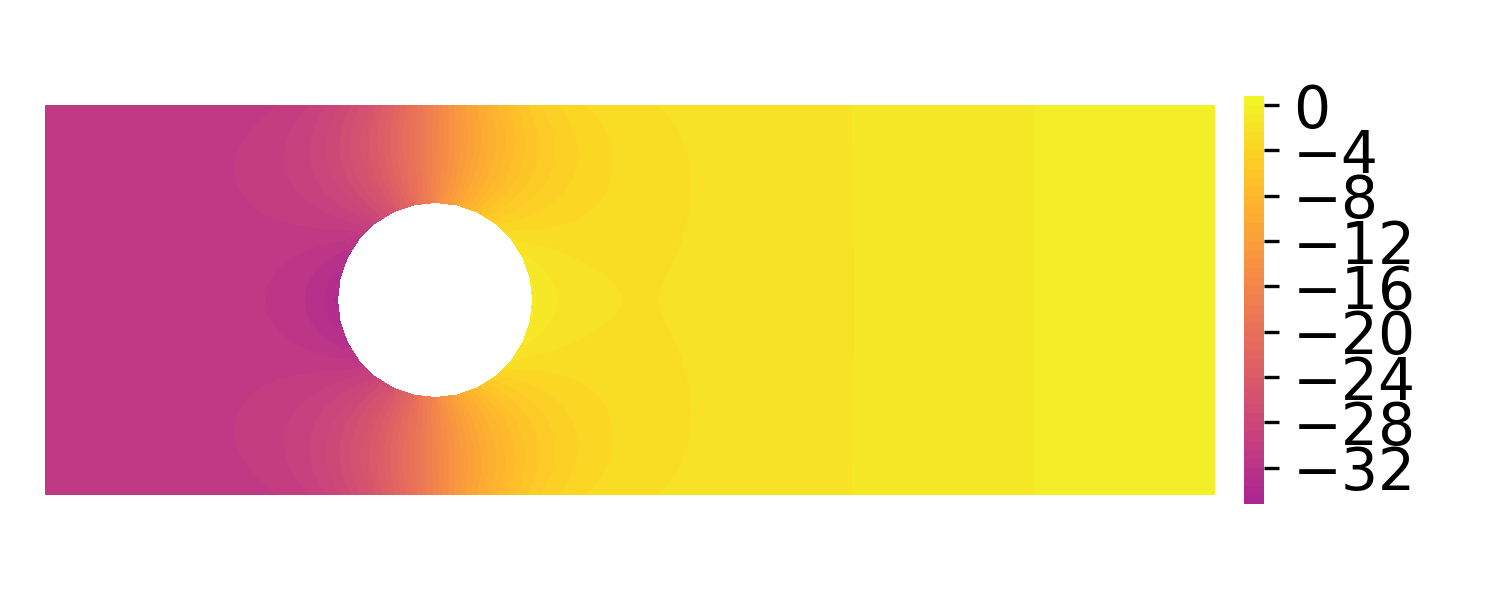} \\[6pt]
            
            \text{AONN} &
            \includegraphics[width=0.43\textwidth]{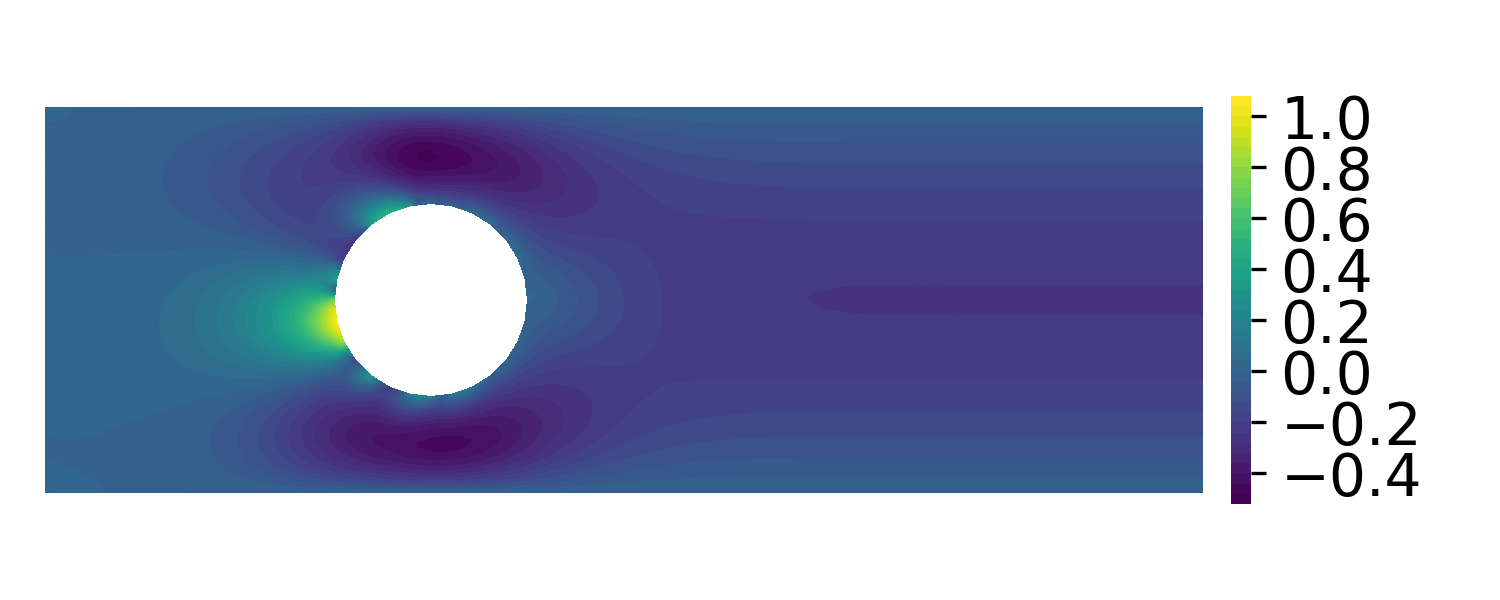} &
            \includegraphics[width=0.43\textwidth]{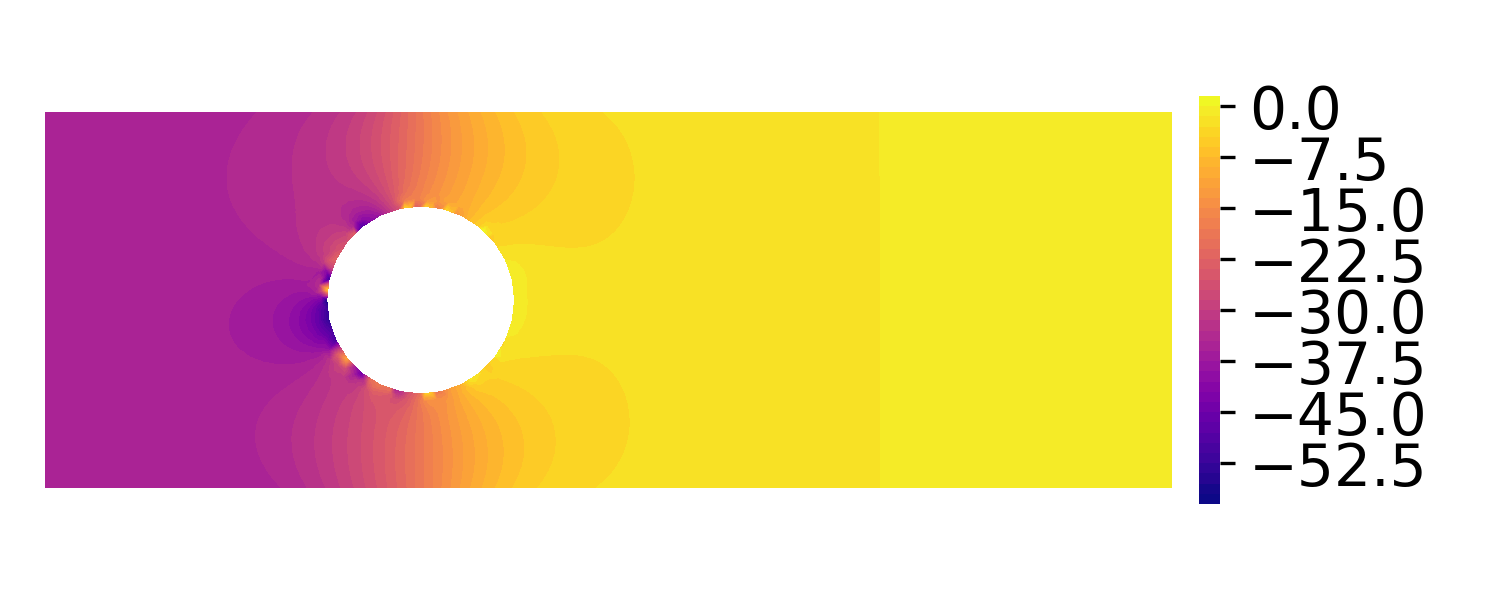} \\[6pt]
            
            \makecell{Adpative\\ AONN} &
            \includegraphics[width=0.43\textwidth]{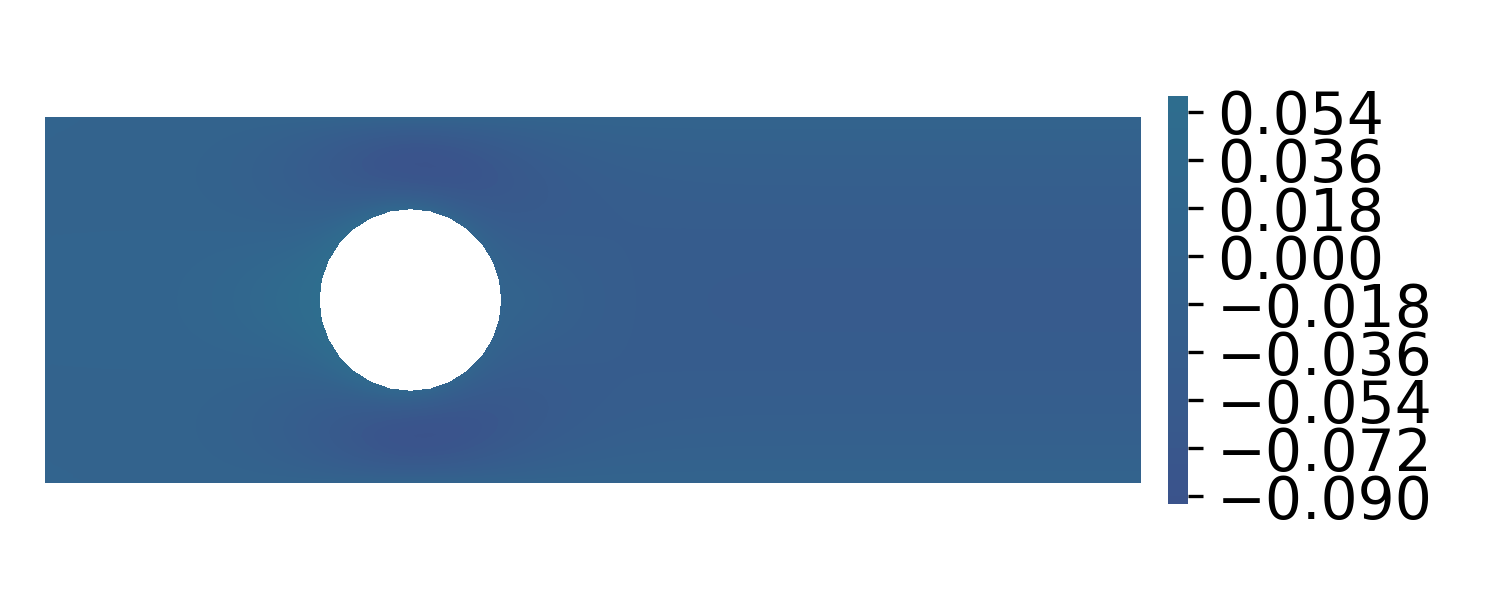} &
            \includegraphics[width=0.43\textwidth]{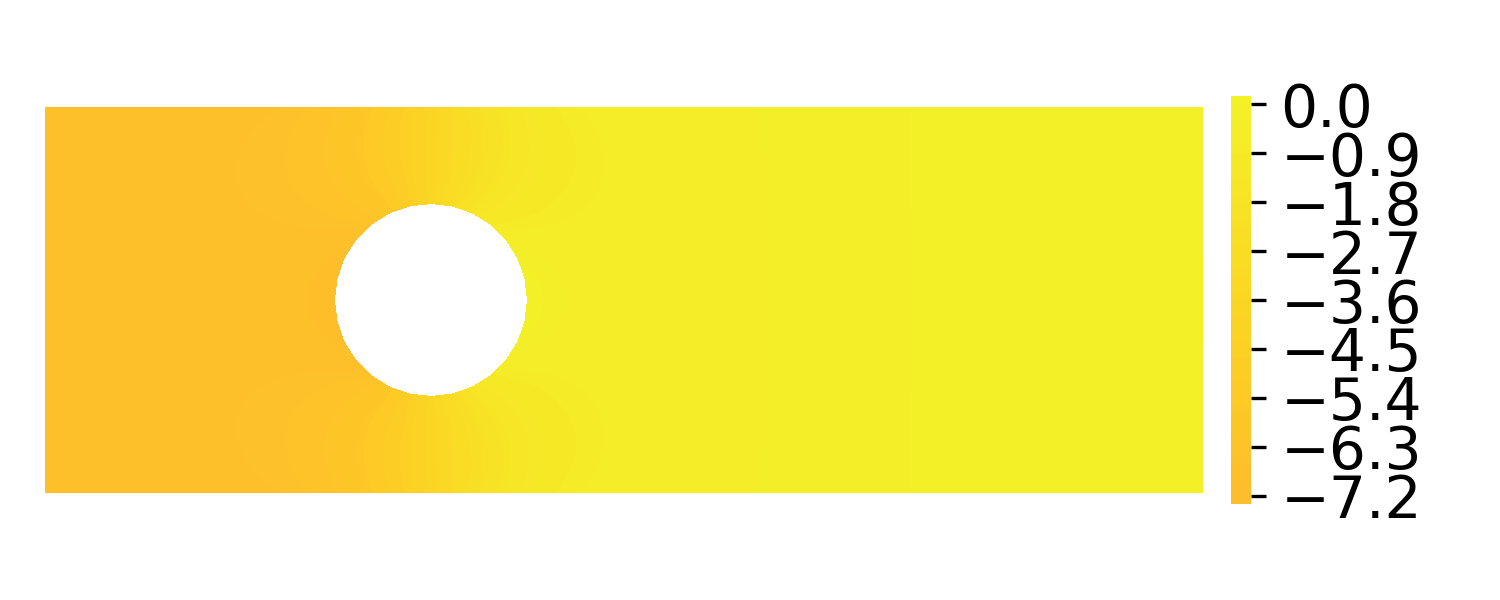}
        \end{tabular}

    \caption{The absolute errors of  velocity field $v(\bx,\xi)$ and  pressure distribution $p(\bx,\xi)$ w.r.t. reference solution for $\xi=10$.}
    \label{stokes_solb}
\end{figure}

Figure~\ref{stokes_solb} displays the absolute errors of $\text{DAS}^2$, AONN, and adaptive AONN with respect to the reference solution. It is clear that adaptive AONN achieves substantially higher accuracy than both $\text{DAS}^2$ and AONN in computing the optimal velocity and pressure fields. When $\xi = 10$, the relative $l_2$ errors of $u(\bx, \xi)$ for $\text{DAS}^2$, AONN, and adaptive AONN are $0.387$, $0.445$, and $0.073$, respectively, while the relative $l_2$ errors of $p(\bx, \xi)$ for the three methods are $0.629$, $0.765$, and $0.119$, respectively. As $\xi$ increases from $10$ to $1000$, the relative errors for both the optimal velocity and pressure fields vary little, with less than a $1\%$ change across the entire parameter range.

\subsection{Test 3: Optimal control for the Laplace equation with parameter-oriented singularity} 
In this test problem, we are going to solve the following parametric OCP:
\begin{equation*}
    \min_{y(r,\theta,\bxi), u(r,\theta,\bxi)} J\big(y(r,\theta,\bxi), u(r,\theta,\bxi)\big)= \frac{1}{2}\left\|y(r,\theta,\bxi)-y_d(r,\theta,\bxi)\right\|_{2, \Omega}^2+\frac{\alpha}{2}\|u(r,\theta,\bxi)\|_{2, \partial \Omega}^2,
\end{equation*}
subject to the following Laplace equations:
\begin{equation*}
\left\{
    \begin{aligned}
        -\Delta y(r,\theta,\bxi)&=0 &&\text{in} \,\Omega, \\
        y(r,\theta,\bxi)&=u(\theta,\bxi) + f(r,\theta,\bxi)  &&\text{on} \,\partial \Omega,
    \end{aligned}
\right.
\end{equation*}
where $\alpha = 0.01$, $\Omega= \mathbb{B}\big((0,0),1\big)$ is a unit circle, and $\bxi$ is a $10$-dimensional parameter with each dimension ranges from $-1$ to $1$. $u(\theta,\bxi)$ denotes the control variable applied on the boundary of the physical domain. The joint spatioparametric space considered in this example is defined as

\begin{equation*}
    \Omega_{\Gamma} \coloneqq \{ (r,\theta, \xi_i) | 0 \leq r \leq 1, 0 \leq \theta \leq 2\pi, -1 \leq \xi_i \leq 1\}  \quad i = 1,2,\cdots,10. 
\end{equation*}
In addition, we impose the following constraint on $u(\theta, \bxi)$:
\begin{equation*}
    0 \leq u(\theta, \bxi) \leq 1 \quad \text { a.e. on } \partial\Omega.
\end{equation*}

The total derivative of $J$ with respect to $u$ is $\mathrm{d}_uJ = \alpha u + \dfrac{\partial p}{\partial \boldsymbol{n}}\Big|_{\partial \Omega}$, where $p$ is the
solution of the corresponding adjoint equation:
\begin{equation*}
\left\{
    \begin{aligned}
        - \Delta p(r,\theta,\bxi) &= y(r,\theta,\bxi) - y_d(r,\theta,\bxi)  \qquad &&\text{in} \, \Omega, \\
        p(r,\theta,\bxi) &= 0 \qquad &&\text{on} \, \partial \Omega. \\
    \end{aligned} 
\right.
\end{equation*} 
If the desired state $y_d(r,\theta,\bxi)$ and the source term $f(r,\theta,\bxi)$ are defined as
\begin{equation*}
    \begin{aligned}
        y_d(r,\theta,\bxi) &= (\frac{1}{2}+\frac{r^2}{2}\cos 2\theta)\,g(\bxi),\\
        f(r,\theta,\bxi) &= \alpha(1+2r^2 \cos 2\theta)\,g(\bxi),
    \end{aligned}
\end{equation*}
where $g(\bxi) = \text{e}^{-10\Vert \bxi \Vert_2^2}$. Then, the exact  solution to this problem is given by
\begin{equation*}
    \begin{aligned}
        y^*(r,\theta,\bxi) &= (\frac{1}{2}+\frac{r^2}{2}\cos 2\theta)\,g(\bxi), \\
        u^*(\theta,\bxi) &= \mathbf{P}_{[0,1]}\big(\cos^2\theta \cdot g(\bxi)\big), \\
        p^*(r,\theta,\bxi) &= \big( -\frac{\alpha}{4}(r^2-1)+(\frac{\alpha}{4}r^2-\frac{\alpha}{4}r^4)\cos 2\theta \big)\,g(\bxi),
    \end{aligned}
\end{equation*}
where $\mathbf{P}_{[u_a,u_b]}$ is the pointwise projection operator onto the interval $[u_a,u_b]$.

To solve this OCP with adaptive AONN, we construct three neural networks: $\hat{y}(r,\theta,\bxi;\theta_y), \hat{u}(r,\theta,\bxi;\theta_u)$ and $\hat{p}(r,\theta,\bxi;\theta_p)$. Specifically, each network consist of $6$ fully connected layers with $20$ neurons in each hidden layer.  For KRnet, we set $K = 3$ and take $L = 6$ affine coupling layers. Each affine coupling layer has $2$ fully connected layers with 64 neurons in each hidden layer. 

The inputs in Algorithm~\ref{alg:DAS-AONN} are set as follows:
the decay factor $\gamma$ is set to $1$; the initial step size $c^{(0)}$ is set to $100$; the initial number of epochs $n^{(0)}$ is set to $200$; the maximum number of epochs $N_{\text{ep}}$ is set to $2000$; the batch size $m$ is set to $4000$; the number of adaptive iterations $N_{\text{adaptive}}$ is set to $10$; and $n_{\text{aug}}$ is set to $0$, implying that a fixed number of training epochs is employed in this testcase. In addition, the three neural networks are trained using the BFGS method, while KRnet is trained using the ADAM optimizer with a learning rate of $0.0001$.

To evaluate the performance of adaptive AONN method, we again conduct a comparative analysis with $\text{DAS}^2$ and AONN. The neural network configurations for $\text{DAS}^2$ and AONN are identical to those used in adaptive AONN.

\begin{figure}[htbp]
\centering
\includegraphics[width=0.9\textwidth]{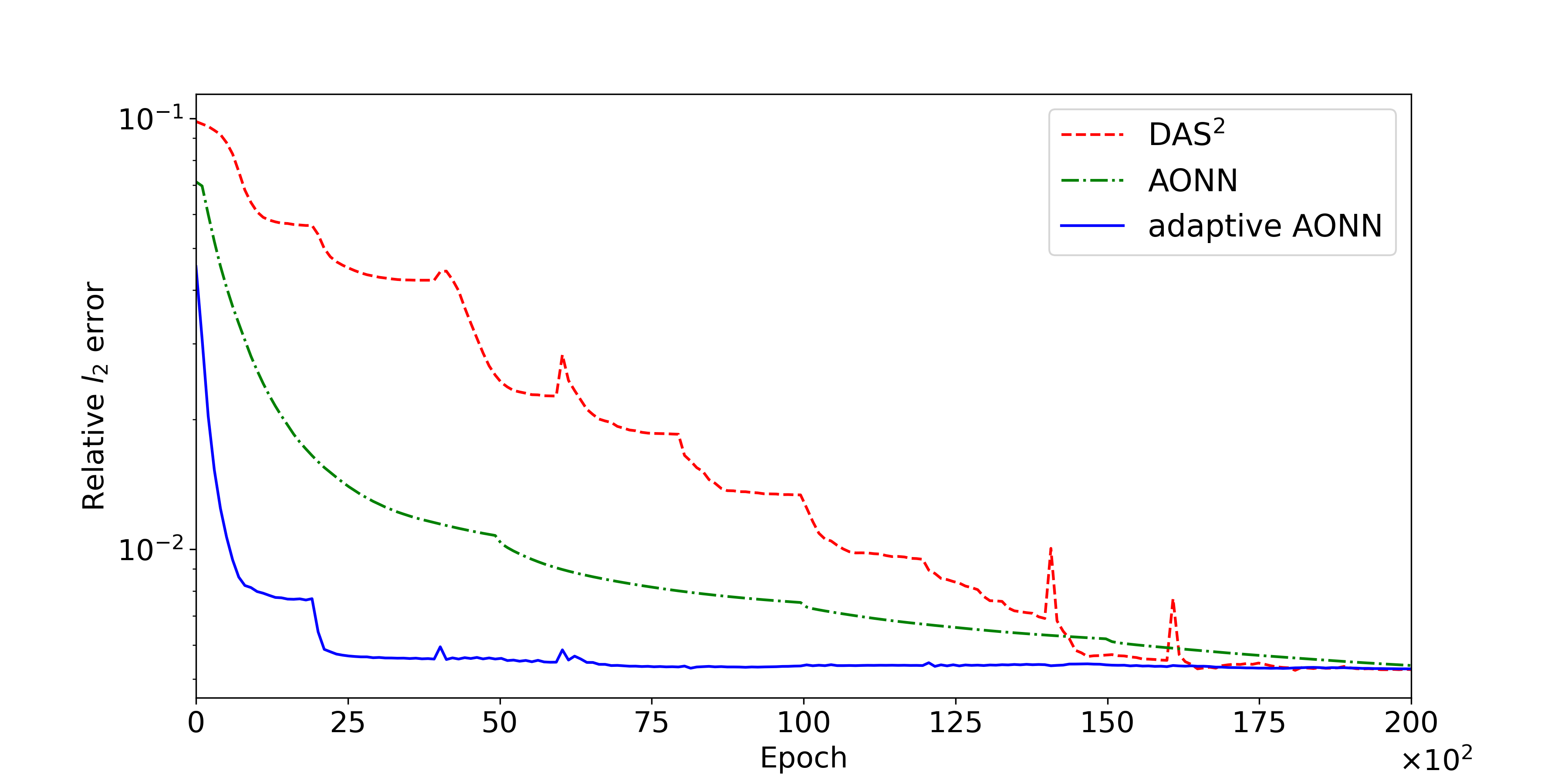}
\caption{Relative $l_2$ error w.r.t. training epochs in Test 3.}\label{test:sigular-compare-l2}
\end{figure}
Figure~\ref{test:sigular-compare-l2} illustrates the evolution of the relative $l_2$ loss for three methods as the number of training epochs increases. Here, the initial training set size is $n_r =|\mathrm{S}_{\Omega,0}| = 6 \times 10^3$. From the figure, it can be seen that adaptive AONN achieves a rapid and stable error reduction from the beginning of the adaptive iterations. While AONN exhibits a monotonically decreasing error trend throughout the training process, yet its convergence is notably slower than that of adaptive AONN. In contrast, the error reduction of $\text{DAS}^2$ appears more irregular. While the overall trend is decreasing, some adaptive stages exhibit pronounced error drops, whereas others show only gradual improvement. Compared with the error curve of $\text{DAS}^2$ in Test 1, the curve of this test case is noticeably smoother. This difference can be explained by the stronger singularity of the solution in Test 3, which makes the effect of adaptive sampling more pronounced.

\begin{figure}[htbp]
    \centering
\subfloat[Relative $L_2$ error w.r.t. sample size.]{
\includegraphics[width=0.45\textwidth]{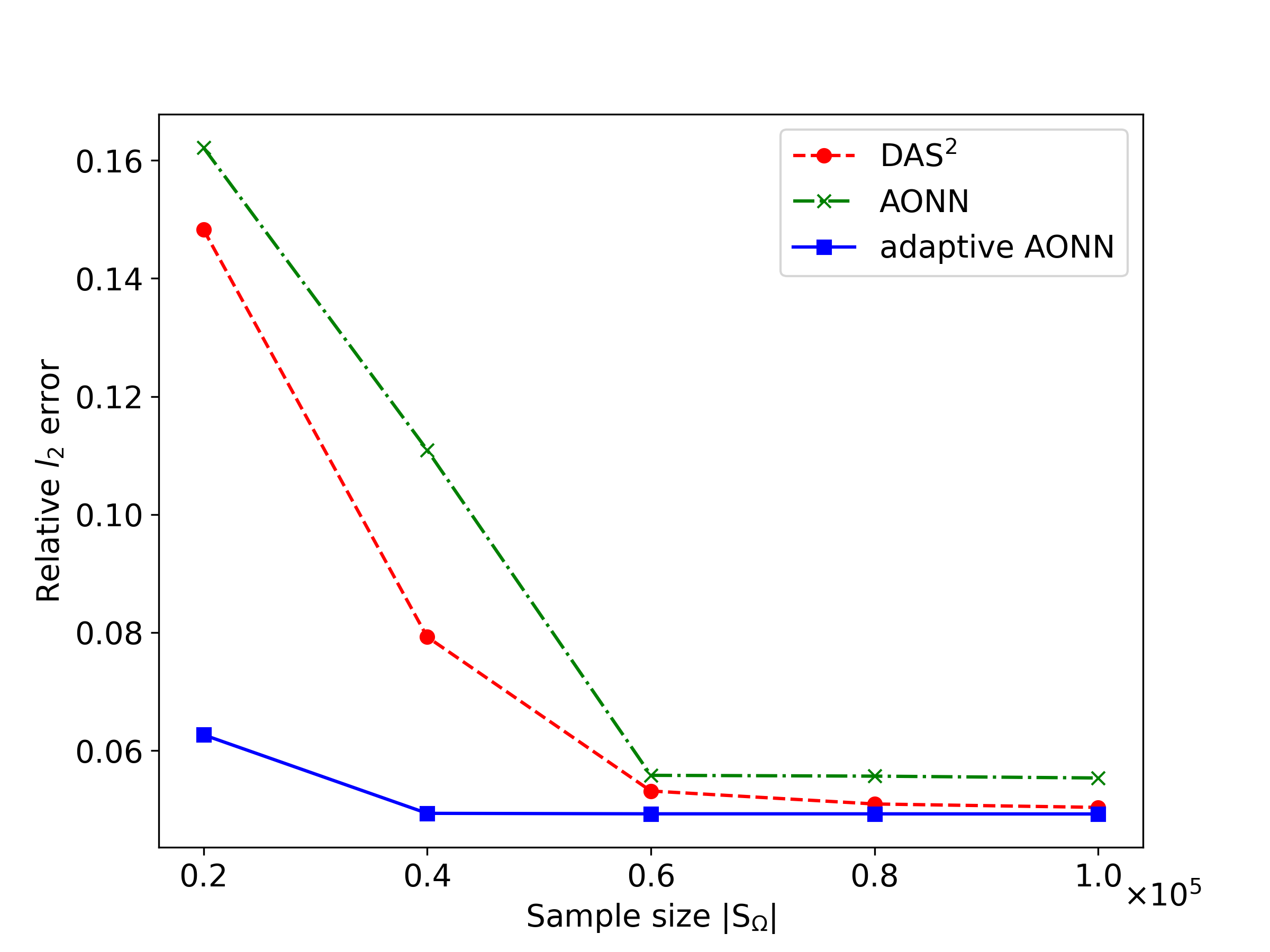}
\label{test:sigular-size-error}
}
\subfloat[Training time w.r.t. sample size.]{
\includegraphics[width=0.45\textwidth]{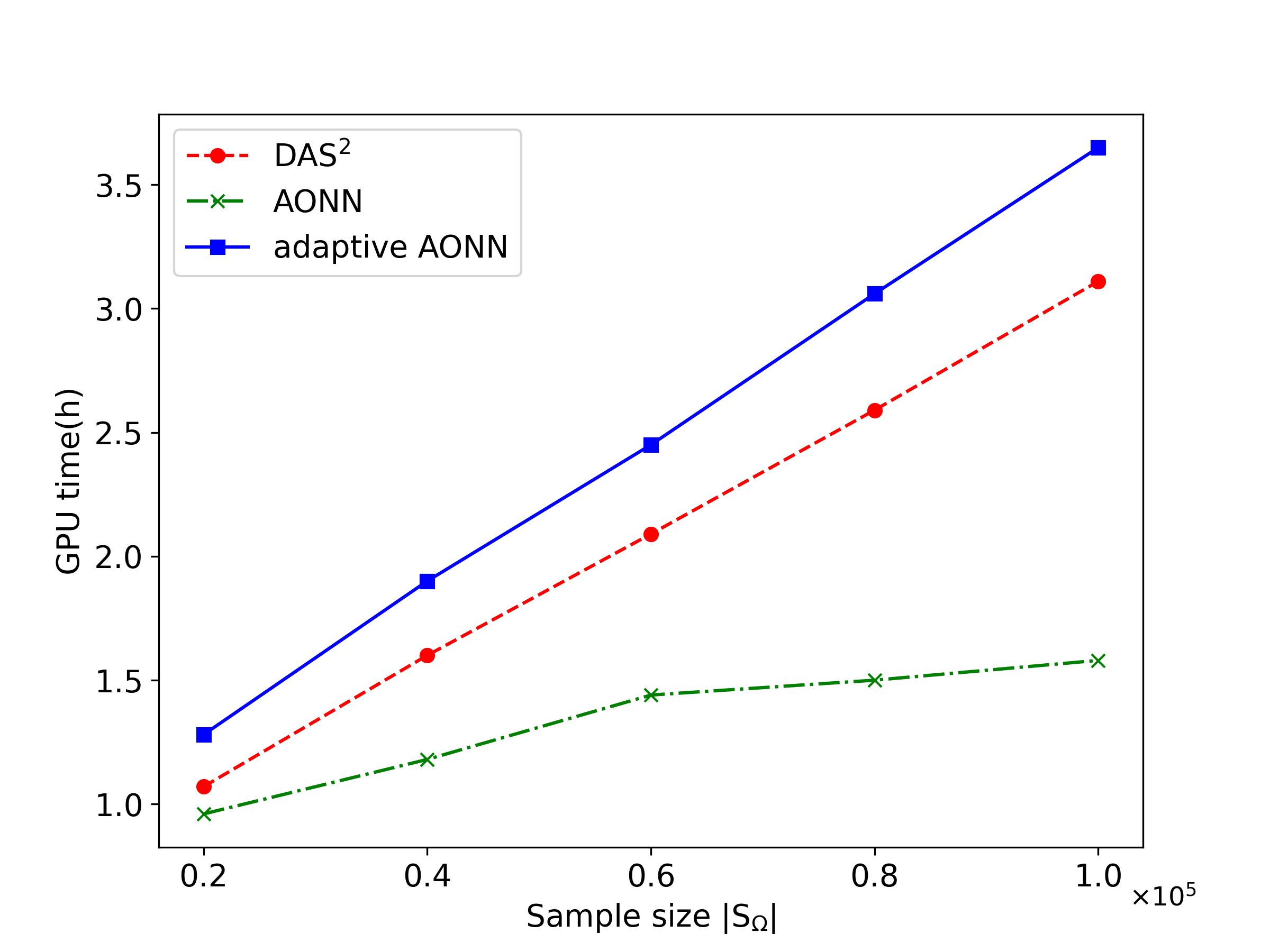}
\label{test:sigular-size-time}
}
\caption{Comparison of $\text{DAS}^2$, AONN and adapive AONN in Test 3.}\label{test:sigular-error}
\end{figure}
We also investigated the relative errors with respect to the sample size and the training time. To do this, we set the initial training set size $\vert \mathrm{S}_{\Omega, 0} \vert = n_r$  to $2\times10^3, 4\times10^3, 6\times10^3, 8\times10^3$, and $1\times10^4$ for adaptive AONN and $\text{DAS}^2$, and resulting in total training size $2\times10^4, 4\times10^4, 6\times10^4, 8\times10^4$, and $1\times10^5$, respectively. The number of collocation points for AONN is set to $10n_r$ directly. The number of adaptivity iterations is $N_{\text{adaptive}} = 10$ with the training number of epochs $N_{\text{ep}} = 3000$.

\begin{table}[!htb]
\centering
\caption{Test 3: Training time and relative error w.r.t. sample size} 
\newcolumntype{C}{>{\centering\arraybackslash}X}
\begin{tabularx}{0.9\linewidth}{CCCCCCC}
  \toprule
  \multirow{2}{*}{$|\mathrm{S}_{\Omega}|$} & \multicolumn{2}{c}{$\text{DAS}^2$} & \multicolumn{2}{c}{AONN} & \multicolumn{2}{c}{adaptive AONN} \\
  \cmidrule(lr){2-3} \cmidrule(lr){4-5} \cmidrule(lr){6-7}
  & Time & Error & Time & Error & Time & Error \\
  \midrule
    $2.0\times10^4$ & 1.07h & 0.14827 & 0.96h & 0.16216 & 1.28h & 0.06268 \\
    $4.0\times10^4$ & 1.60h & 0.07926 & 1.18h & 0.11073 & 1.90h & 0.04937 \\
    $6.0\times10^4$ & 2.09h & 0.05316 & 1.44h & 0.05581 & 2.45h & 0.04929 \\
    $8.0\times10^4$ & 2.59h & 0.05096 & 1.50h & 0.05568 & 3.06h & 0.04929 \\
    $1.0\times10^5$ & 3.11h & 0.05038 & 1.58h & 0.05536 & 3.65h & 0.04927 \\
  \bottomrule
\end{tabularx}
\label{sample_table_para}
\end{table}

The results can be seen in Figure~\ref{test:sigular-error} and Table~\ref{sample_table_para}. 
From Figure~\ref{test:sigular-size-error}, we observe results similar to those of examples in which the singularity arises from the geometric space. For adaptive AONN, the constructed neural networks effectively capture the singularity with only a few collocation points, whereas the other two methods require significantly more samples to achieve convergence. Figure~\ref{test:sigular-size-time} also presents the training time with respect to sample size. The training times of $\text{DAS}^2$ and adaptive AONN are comparable, while AONN has the shortest training time. However, due to the substantial advantage in sample efficiency for achieving the same level of accuracy, adaptive AONN may outperform both AONN and $\text{DAS}^2$ in terms of overall computational cost. This is further confirmed by Table~\ref{sample_table_para}, which shows that adaptive AONN takes only 1.28 hours to reach a comparable accuracy, whereas $\text{DAS}^2$ and AONN require 2.09 and 1.44 hours, respectively.

\section{Conclusions}\label{section_conclude}
In this paper, we have proposed a hybrid methodology, adaptive AONN, for constructing surrogate models to solve optimal control problems constrained by parametric partial differential equations with singular features. The key idea of adaptive AONN is to adaptively update the surrogate models and their training set alternately, with the training set update achieved by approximating the probability distribution induced by the residual profile. At each stage, the solution obtained through AONN enhances the ability of $\text{DAS}^2$ to approximate the target distribution, while the samples generated by $\text{DAS}^2$ aid in training the surrogate model within AONN, creating a mutually reinforcing computational synergy. Adaptive AONN demonstrates strong performance in handling singularity problems, effectively constructing surrogate models that accurately approximate solutions near singular regions.












\bibliographystyle{elsarticle-num} 
\bibliography{AdaptiveAONN_ref}
\end{document}